\documentclass[12pt]{article} 
\usepackage{amsmath, amssymb}
\usepackage{amsthm} 
\usepackage[all]{xy}

\theoremstyle{plain} 

\theoremstyle{definition} 

\def\L{{\cal L}}
\def\I{{\cal I}}
\def\bgn{\begin}

\def\CL{\text{\rm CL}}

\def\J{{\cal J}}
\def\L{{\cal L}}
\def\Ob{\text{\rm Ob}}
\def\1{{[1]}}
\def\2{{[2]}}
\def\3{{[3]}}
\def\({\left(}
\def\){\right)}
\def\s-circ{\,{\scriptstyle{\circ}}\,}
\def\<<{<\negthinspace \negthinspace<}
\def\Ad{\text{\rm Ad}}
\def\ad{\text{\rm ad}}

\def\bgn{\begin}

\def\endaln{\end{align}}

\def\bE{\bold E}
\def\<{<\negthinspace \negthinspace <}

\def\({\left(}
\def\){\right)}

\def\[{\big[\neg\big[}
\def\]{\big]\neg\big]}
\def\al{\al}

\def\M{{\cal M}}

\def\a{\alpha}
\def\b{\beta}

\def\e{\varepsilon}
\def\gam{\gamma}
\def\Gam{\Gamma}
\def\k{\kappa}
\def\del{\delta}
\def\lam{\lambda}
\def\ome{\omega}
\def\Ome{\Omega}
\def\sig{\sigma}
\def\Sig{\Sigma}

\def\Q{\Bbb Q}
\def\R{\Bbb R}
\def\C{\Bbb C}

\def\O{\Bbb O}
\def\P{\frak P}
\def\M{\frak M}

\def\w{\wedge}
\def\({\left(}
\def\){\right)}

\def\neg{\negthinspace}
\def\h{\hat}

\def\til{\tilde}
\def\wtil{\widetilde}
\def\ch{\check}

\def\ol{\overline}
\def\pa{\partial}

\def\ran{\rangle} 
\def\lan{\langle}

\def\ss{\scriptscriptstyle}

\def\arrow{\longrightarrow}

\def\SL{\ss{\text{\rm SL}}}

\def\bsh{\backslash}

\def\:{\, :\,}
\def\CL{\text{\rm CL}}
\def\TT{T\oplus T^*}
\def\complex{generalized complex }
\def\K\"ahler{generalized K\"ahler}

\def\10{\displaystyle L^{10}}
\def\2{\displaystyle L^2}
\def\c0{\displaystyle C^0}

\def\10{\displaystyle L^{10}}
\def\2{\displaystyle L^2}
\def\del{\delta}
\def\del2{\displaystyle L^2_{0,\delta}}
\def\c0{\displaystyle C^0}

\def\del{\delta}

\def\K{{\cal K}}
\def\M-A{\text{\rm Monge-Amp\`ere}}
\def\O{{\cal O}}

\def\M-A{\text{\rm Monge-Amp\`ere}}
\def\[{\big[\,}
\def\]{\,\big]}
\def\End{\text{\rm End\,}}

\def\id{\text{\rm id}}
\def\P{\Bbb P}
\title{{Unobstructed K-deformations of Generalized Complex Structures and Bihermitian Structures}} 
\author{
\textsc{Ryushi Goto$^{*}$} 
}
\date{} 
\begin{document}
\maketitle
\footnote{ 
2010 \textit{Mathematics Subject Classification}.
Primary 53D18; Secondary 53C55.
}
\footnote{ 
\textit{Key words and phrases}. 
generalized complex manifolds, generalized K\"ahler structures, bihermitian structures, unobstructed deformations}
\footnote{ 
$^{*}$Partly supported by the Grant-in-Aid for Scientific Research (C),
Japan Society for the Promotion of Science. 
}
\begin{abstract}
We introduce ``K-deformations" of generalized complex structures on a compact K\"ahler manifold $M=(X, J)$ with an effective anti-canonical divisor and 
show that obstructions to K-deformations of generalized complex structures on $M$ always vanish.
Applying unobstructed K-deformations and the stability theorem of generalized K\"ahler structures, 
 we construct deformations of bihermitian structures in the form 
$(J, J^-_t, h_t)$ on a compact K\"ahler surface with a non-zero holomorphic Poisson structure.
Then we prove that a compact K\"ahler surface
$S$ admits a non-trivial bihermitian structure with the torsion condition and the same orientation
if and only if $S$ has a non-zero holomorphic Poisson structure. 
Further we obtain bihermitian structures $(J, J^-, h)$ on del Pezzo surfaces, degenerate del Pezzo surfaces and 
some ruled surfaces for which the complex structure $J$ is not equivalent to $J^-$ under diffeomorphisms. 
\end{abstract}
\tableofcontents
\numberwithin{equation}{section}
\section{Introduction}
The concept of generalized complex structures were introduced by Nijel Hitchin \cite{Hi1}
depending on a simple idea by replacing the tangent bundle on a manifold with the direct sum of 
the tangent bundle and the cotangent bundle which include
both symplectic and complex structures as special cases. Interesting generalized complex structures 
are arising as hybrid of symplectic and complex structures.
An associated notion of generalized K\"ahler structures consists of two commuting generalized complex structures
which yield a generalized metric. It is remarkable that generalized K\"ahler structures are equivalent to the so-called bihermitian structures \cite{Gu1}.
Bihermitian structures on complex surfaces are extensively studied 
from view points of conformal geometry and complex geometry \cite{A.D, A.G.G, F.P}.  Beihermitian structures are also appeared in the Physics as geometric structures on the target space for the supersymmetric $\sigma$-model \cite{Ro}.  
Main purpose of this paper is to solve remarkable problems in bihermitian geometry by using
deformations of generalized K\"ahler structures.
\\
Let $M=(X,J)$ be a compact K\"ahler manifold with effective anti-canonical divisor where $X$ is the underlying differential manifold. 
Then there is the \complex structure $\J_J$ induced from the complex structure $J$.
At first we shall show that there exist  
certain unobstructed deformations of \complex structures of $\J_J$, which is called
{\it $K$-deformations}.
This is regarded as a generalization of unobstructedness theorem of Calabi-Yau manifolds \cite{Bo}, \cite{Ti}
(see also Miyajima's result on
unobstructed deformations in the case of normal isolated singularity \cite{Miya}.)
We apply a unified method as in \cite{Go-1}, \cite{Go0} to a meromorphic $n$-form with a pole along the anti-canonical divisor and show that the obstructions to deformations vanish:\\ \\
\noindent
{\bf Theorem} \ref{th: unobstructedness  theorem-1} {\it
Let $M=(X,J)$ be a compact K\"ahler manifold of dimension $n$
and we denote by $\J_J$ the \complex structure given by $J$.
If $M$ has an effective, anti-canonical divisor $D$, 
then $M$ admits unobstructed $K$-deformations of \complex structures 
$\{ \J_t\}$ starting with $\J_0=\J$ which are 
parametrized by an open set of $H^2(\ol L(-D))\cong H^{n,2}\oplus H^{n-1,1}\oplus H^{n-2, 0}$.}
\\ 
In the cases of K\"ahler surfaces, we obtain the following: \\ \\
{\bf Corollary }\ref{cor: unobstructedness S}  {\it 
Let $S$ be a compact K\"ahler surface with the complex structure $J$ and 
a K\"ahler form $\ome$. 
If $S$ has an effective, anti-canonical divisor $[D]=-K_S$, 
then $S$ admits unobstructed $K$-deformations of \complex structures parametrized by an open set of
the full cohomology group $H^0(S)\oplus H^2(S)\oplus H^4(S)$ of even degree.}
\\  \\
The obstruction space of K-deformations on compact K\"ahler surface $S$ is given by $H^{1,2}\oplus H^{0,1}\cong H^1(S, \C)$. 
Thus the obstruction space does not vanish if $b_1(S)\neq 0$. 
For instance, the product of $\C P^1$ and the elliptic curve $E$ admits unobstructed K-deformations, nevertheless
the obstruction space does not vanish.

We apply our unobstructed K-deformations to construct bihermitian structures on compact K\"ahler surfaces.
A bihermitian structure on a differential manifold $X$ is a triple $(J^+, J^-, h)$ consisting of two complex structure $J^+$ and $J^-$ and 
a metric $h$ which is a Hermitian metric with respect to both $J^+$ and $J^-$. 
In this paper we always assume that a bihermitian structure satisfies the torsion condition: 
\bgn{equation}\label{eq: torsion condition}
-d^c_+\ome_+ =d^c_-\ome_- =db,
\end{equation}
where $d^c_\pm=\sqrt{-1}(\ol\pa_\pm-\pa_\pm)$ and 
$\ome_\pm$ denote the fundamental 2-forms with respect to $J^\pm$ and 
$b$ is a real $2$-form.
A bihermitian structure on $X$ is {\it non-trivial} if there is a point $x\in X$ such that 
$J^+_x\neq \pm J^-_x$.
A bihermitian structure with $J^+_x\neq J^-_x$ for all $x\in X$ is called 
{\it a strongly bihermitian structure}.
If $J^+$ and $J^-$ induce the same orientation then $(J^+, J^-, h)$ is 
{\it a bihermitian structure with the same orientation.}
A bihermitian structure $(J^+, J^-, h)$ on $X$ is {\it distinct} if 
the complex manifold $(X, J^+)$ is not biholomorphic to the complex manifold $(X, J^-)$. 
We say {\it a complex manifold $M=(X, J)$ admits a bihermitian structure} if there is a bihermitian structure $(J^+, J^-, h)$ on $X$ with $J^+=J$. 
We treat the following question in this paper: 
 $$
   \text{\rm     Which compact complex surfaces admit nontrivial bihermitian structures ?}
 $$
The question was addressed by Apostolov, Gauduchon and Grantcharov \cite{A.G.G}. Kobak \cite{Kob} gave strongly bihermitian structures on the torus $T^4$.
Hitchin constructed (non-strongly) bihermitian structures on del Pezzo surfaces by using the Hamiltonian diffeomorphisms \cite {Hi3} .
Fujiki and Pontecorvo constructed anti-self-dual bihermitian structures on certain non-K\"ahler surfaces in class VII , especially hyperbolic and parabolic Inoue surfaces by using the Twistor spaces \cite{F.P}.
There is a one to one correspondence between bihermitian structures with the torsion condition (\ref{eq: torsion condition})
and generalized K\"ahler structures \cite{Gu1}. 
Thus we can obtain bihermitian structures by constructing generalized K\"ahler structures.
Lin and Tolman developed the generalized K\"ahler quotient construction to obtain examples of generalized K\"ahler manifolds \cite{L.T}
and the author established the stability theorem of generalized K\"ahler structures to construct generalized K\"ahler deformations on a compact K\"ahler manifold with a holomorphic Poisson structure \cite{Go1}, \cite{Go2}.  However these generalized K\"ahler structures do not give a precise answer of the above question since both complex structures of the corresponding bihermitian structures may be deformed.  

If we try to obtain deformations of bihermitian structures 
$(J, J^-_t, h_t)$ fixing one of complex structures, then we encounter a problem of the obstructions to deformations of \complex structures (see \cite{Go3}). 
We use K-deformations instead of ones of \complex structures in order to overcome the difficulty.
Since K-deformations on compact K\"ahler surfaces are unobstructed, 
we obtain 
\medskip \\
\noindent{\bf Theorem }\ref{th: bihermitian on compact Kahler surfaces}    {\it 
Let $S=(X, J)$ be a compact K\"ahler surface with a K\"ahler form $\ome$.
If $S$ has a non-zero holomorphic Poisson structure, then $S$ admits  deformations of non-trivial
bihermitian structures 
$(J, J^-_t, h_t)$ with the torsion condition which satisfies 
$$
\frac{d}{dt}J^-_t|_{t=0}=-2(\b+\ol\b)\cdot\ome
$$
and $J^-_0=J$ and $h_0$ is the K\"ahler metric of $(X, J, \ome)$, where 
$\b\cdot\ome$ is the $T^{1,0}_J$-valued $\ol\pa$ closed form of type $(0,1)$ which gives 
the Kodaira-Spencer class $[\b\cdot\ome]\in H^1(S, \Theta)$ of the deformations $\{J^-_t\}$.
} \\ 
It is shown that a non-trivial bihermitian structure with the torsion condition and the same orientation on a compact surface gives a nonzero holomorphic Poisson structure 
\cite{A.G.G}, \cite{Hi3}, (see Proposition 2 and Remark 2 in \cite{A.D}).
Thus it follows from our theorem \ref{th: bihermitian on compact Kahler surfaces}  that \smallskip\\
{\bf Theorem }\ref{th: bihermitian and Poisson}  {\it  A compact K\"ahler surface admits non-trivial bihermitian structure with the torsion condition and the same orientation if and only if $S$ has nonzero holomorphic Poisson structure.} \smallskip \\
For instance it turns out that all degenerate del Pezzo surface and all Hirzebruch surfaces admit non-trivial bihermitian structures with the torsion condition and the same orientation. 
Further since there is a classification of Poisson surfaces \cite{B.M}, \cite{Sa}, \cite{Sa1},  we obtain all compact K\"ahler surfaces which admit bihermitian structures with the torsion condition and the same orientation.
Let $T^*\Sig$ be the cotangent bundle for every Riemannian surface $\Sig$ with genus $g$ and $S$ the projective space bundle $\P(T^*\Sig\oplus {\cal O}_\Sig)$ with the fibre $\P^1$.
Then $S$ has an effective divisor 
$2[E_\infty]$, where $E_\infty$ is the section of $S\to \Sig$ with intersection number
$2-2g$. 
Let $\b$ be a holomorphic Poisson structure with the zero locus $2[E_\infty]$. 
Then it turns out that the class $[\b\cdot\ome]\in H^1(S, \Theta)$ does not vanish. 
We denote by $X$ the underlying differential manifold of the complex surface $S$.
Thus we obtain \smallskip \\
{\bf Theorem }\ref{th: distinct bihermitian}
{\it There is a family of distinct bihermitian structures $(J, J^-_t, h_t)$ with the torsion condition and the same orientation on 
$S:=\P(T^*\Sig\oplus {\cal O}_\Sig)$, 
that is, the complex manifold $(X, J^-_t)$ is not biholomorphic to $S=(X, J)$ for small $t\neq 0$.}
\smallskip\\ 
In section 2 we obtain unobstructed K-deformations of \complex structures.
In section 3 unobstructed K-deformations of \complex structures
is given by the action of CL$^2(-D)$ which is necessary for our construction of generalized K\"ahler structures.
In section 4 we recall the stability theorem of generalized K\"ahler structures and 
in section 5 we describe a family of sections $\Gam^\pm(a(t),b(t))$ of GL$(TX)$ which gives deformations of 
bihermitian structures $(J^+_t, J^-_t, h_t)$. 
In section 6,7 we shall show our main theorem \ref{th: bihermitian on compact Kahler surfaces}. 

The author would like to thank for Professor A. Fujiki for valuable comments. 
He is thankful to Professor V. Apostolov for meaningful comments on bihermitian structures at Sugadaira conference and 
at Montr\'eal.
He is grateful to Professor M. Gualtieri for a remarkable discussion at Montr\'eal.

\section{Unobstructed K-deformations of \complex structures}
\noindent
Let $M=(X, J)$ be a compact K\"ahler manifold of dimension $n$, where 
$X$ is the underlying differential manifold and $J$ is the complex structure on it.
Then there is the \complex structure $\J$ defined by $J$. We denote by $T^{1,0}$ the holomorphic tangent bundle and $\w^{0,1}$ is the $C^\infty$ vector bundle 
of forms of type $(0,1)$. 
The vector bundle $\ol L$ is defined to be the direct sum $T^{1,0}\oplus \w^{0,1}$ and we denote by
$\w^r\ol L$ the bundle of anti-symmetric tensors of $\ol L$ with degree $r$ which has the decomposition, 

$$
\w^r\ol L= \oplus_{p+q=r} T^{p,0}\otimes\w^{0,q},
$$
where $T^{p,0}$ is  the bundle of $p$-vectors  and $\w^{0,q}$ is  the bundle of forms of type $(0,q)$.
Thus we have the complex $(\w^\bullet\ol L, \ol\pa)$ which is the direct sum of the ordinary 
$\ol\pa$-complexes $(T^{p,0}\otimes \w^{0,\bullet}, \ol\pa)$. 
The sheaf of smooth sections of $\w^r\ol L$ is denoted by $\w^r\ol\L$. 
We assume that $M$ has an anti-canonical divisor $D$ which is given by the zero locus of a holomorphic section $\b\in H^0(M, K^{-1}_M)$. Thus the canonical line bundle $K_M$ is the dual bundle $-[D]$.
We denote by $I_D$ the ideal sheaf of the divisor $D$ which is the sheaf of sections of the canonical line bundle $K_M$. 
We define a sheaf $\w^r\ol\L(-D)$ by 
\bgn{equation}
\w^r\ol\L(-D)(U):=\{\, fa\, |\, f\in I_D(U), \,\, a\in \w^r\ol L(U)\,\,\},
\end{equation}
where $U$ is an open set of the manifold $X$.
In particular, the sheaf $\w^0\ol\L(-D)$ is given by $I_D\otimes C^\infty(X)$. 
We denote by $\Gam(X, \w^r\ol \L(-D))$ the set of global smooth sections of the sheaf $\w^r\ol\L(-D)$.
(In this paper, for simplicity, we often say a section of a sheaf instead of a global section of a sheaf.)
The sheaf $\w^r\ol\L(-D)$ is locally free which is the sheaf of smooth sections of a vector bundle 
$\w^r\ol L(-D)$. Tensoring with $\b^{-1}$ gives an identification, 
$$
\w^r\ol L(-D)\cong \w^r \ol L\otimes[-D]\cong \w^r\ol L\otimes K_M
$$
Then we have 
\bgn{align*}
\w^0  \ol L(- D)=&\w^{0,0}\otimes K_M\\
\w^1\ol L(-D)=&(\w^{0,1}\otimes K_M)\oplus ( T^{1,0}\otimes K_M)\\
\w^2\ol L(-D)=&(\w^{0,2}\otimes K_M)\oplus 
(T^{1,0}\otimes \w^{0,1}\otimes K_M)\oplus (T^{2,0}\otimes K_M)
\end{align*}
Thus we have the subcomplex $(\w^\bullet\ol \L(-D), \ol\pa)$ of the complex $(\w^\bullet\ol \L, \ol\pa)$.
Let $\Ome$ be a meromorphic $n$-form with simple pole along the divisor $D$ which is unique up to 
constant multiplication since $H^0(M, K\otimes[D])\cong \C$.
Note that the $d$-closed form $\Ome$ is locally written as
$$
\Ome|_U = \frac{dz_1 \w\cdots\w dz_n}{f},
$$
by a system of holomorphic coordinates $(z_1, \cdots, z_n)$ on a small open set $U$ and 
$f\in I_D(U)$.
The action of $\w^r\ol \L(-D)$ on $\Ome$ is defined by the interior and exterior product which is the spin representation of the Clifford algebra $\CL(\TT)$.
Then it turns out that this action yields a smooth differential form on $M$. 
Thus we obtain an identification, 
$$
\w^r\ol L(-D)\cong \oplus_{p+q=r}\w^{n-p, q},
$$
where $\w^{n-p,q}$ is the $C^\infty$ vector bundle of forms of type $(n-p, q)$.
\bgn{lemma} The direct sum $\oplus_r\w^r \ol \L(-D)$ is involutive with respect to the Schouten bracket.
\end{lemma}
\bgn{proof}
A local section $\a$ of the sheaf  $\w^r\ol \L$ is a local section of the sheaf $\w^r\ol \L(-D)$ if and only if $\a\cdot\Ome $ is a smooth differential form. 
The Schouten bracket$[\a_1,\a_2]_S$ for $\a_1, \a_2\in \w^\bullet\ol\L(-D)$ is given by 
$[[d,\a_1]_G, \a_2]_G$, where $[\, ,\,]_G$ denotes the graded bracket, which is the derived bracket construction
(see \cite{Go1}, \cite{Kos}).
Thus $[\a_1, \a_2]_S\cdot\Ome$ is a smooth differential form since $d\Ome=0$. 
It follows that $[\a_1, \a_2]_S$ is a local section of $\w^r\ol \L(-D)$.
\end{proof}

We define a vector bundle $U^{-n+r}$ by 
$$
U^{-n+r}:=\bigoplus_{\stackrel {\scriptstyle  p+q=r}{0\leq p, q\leq n}}\w^{n-p, q}
$$
Then the identification $\w^r\ol L(-D)\cong \oplus_{p+q=r}\w^{n-p, q}$ gives an isomorphism between complexes, 
$$
(\w^\bullet L(-D), \ol\pa)\cong (U^{-n+\bullet}, \ol\pa).
$$
Thus the cohomology group $H^\bullet(\ol L(-D))$ is given by the direct sum of the Dolbeault cohomology groups,
$$
H^r(\ol L(-D))\cong  \bigoplus_{\stackrel {\scriptstyle  p+q=r}{0\leq p, q\leq n}}H^{n-p, q},
$$
where $H^{n-p, q}$ is the Dolbeault cohomology group $H^{n-p, q}(M)$.
Thus we have
\bgn{proposition}\label{prop: H(ol L(-D))}
\bgn{align*}
H^0(\ol L(-D) )=&H^0(M, K_M)\cong H^{n,0}\\ 
H^1(\ol L(-D)) =&H^1(M,K_M)\oplus H^0(M,\Theta\otimes K_M)\cong H^{n,1}
\oplus H^{n-1,0}\\ 
H^2(\ol L(-D)) =&H^2(M, K_M)\oplus H^1(M, \Theta\otimes K_M)\oplus H^0(M,\w^2\Theta\otimes K_M)\\
\cong& H^{n,2}\oplus H^{n-1,1}\oplus H^{n-2, 0}
\end{align*}
\end{proposition}
We define  vector bundles $\bE^\bullet$ by 
\bgn{align*}
&\bE^{-1}=U^{-n},\qquad\quad \,\,  \bE^0=U^{-n+1}, \\
& \bE^1=U^{-n}\oplus U^{-n+2}, 
\bE^2=U^{-n+1}\oplus U^{-n+3},\, 
\end{align*}
Then we have the complex by using the exterior derivative $d$,
$$
0\arrow \bE^{-1}\overset{d}\arrow  \bE^0\overset{d}\arrow \bE^1\overset{d}\arrow \bE^2\overset{d}\arrow\cdots
$$
The cohomology groups $H^\bullet(\bE^\bullet)$ of the complex $(\bE^\bullet, d)$ is given by the direct sum of the Dolbeault cohomology groups,
\bgn{align*}
&H^{-1}(\bE^\bullet)= H^{n,0}\\
&H^0(\bE^\bullet) = H^{n,1}\oplus H^{n-1,0}\\
&H^1(\bE^\bullet)= H^{n,2}\oplus H^{n-1,1}\oplus H^{n-2, 0}\oplus H^{n,0}\\
&H^2(\bE^\bullet)=H^{n,3}\oplus H^{n-1,2}\oplus H^{n-2, 1}\oplus H^{n-3,0}\oplus H^{n,1}\oplus H^{n-1,0}
\end{align*}

Let $\e_1$ be a smooth global section of $\w^2\ol\L(-D)\oplus\w^0\ol \L(-D)$ with 
$d\e_1\cdot\Ome=0$. 
For such a section $\e_1$ of $\w^2\ol\L(-D)\oplus \w^0\ol\L(-D)$,
we shall construct a family of smooth global sections $\e(t)$ of $\w^0\ol \L(-D)\oplus \w^2\ol \L(-D)$ which gives 
deformations of maximal isotropic subbundles $\{\ol L_t\}$ by the Adjoint action
\bgn{equation}\label{eq: Lt}
L_t:=\Ad_{e^{\e(t)}}L=\{\, E+[\e(t), E]\, |\, E\in L\,\}.
\end{equation}
Then the decomposition $(\TT)^\C=L_t\oplus \ol L_t$ gives an almost \complex structure $\J_t$ whose eigenspaces are $L_t$ and $\ol{L_t}$, where note that $\ol L_t$ is the complex conjugate of $L_t$.
Thus a family of section $\{\e(t)\}$ yields deformations of almost \complex structures $\{\J_t\}$,
where $t$ is a parameter of the deformation.
A family of section $\e(t)$
is given in the form of power series,
$$
\til\e(t)=\e_1t+\e_2\frac{t^2}{2!}+\e_3\frac{t^3}{3!}+\cdots.
$$
If $\J_t$ is integrable, deformations $\J_{t}$ given by  a family of smooth global sections $\e(t)$ of $\w^0\ol \L(-D)\oplus \w^2\ol \L(-D)$ is called {\it K-deformations of \complex structures}.
The structures $\J_{t}$ are integrable if and only if  the family of global sections $\e(t)$ satisfies the Maurer-Cartan equation, 
$$
\ol\pa\e(t)+\frac12[\e(t), \e(t)]_S=0,
$$
where $[\e(t), \e(t)]_S\in \Gam (X, \w^3\ol \L(-D))$ denotes the Schouten bracket of $\e(t)$.
The action of $e^{\e(t)}$ on $\Ome$ gives a non-degenerate, pure spinor 
$e^{\e(t)}\cdot\Ome$ which induces the almost \complex structure $\J_t$.
It is crucial to solve the following equation, 
\bgn{equation}\label{eq: de e(t) Ome=0}
de^{\e(t)}\cdot\Ome=0,
\end{equation}
rather that the Maurer-Cartan equation. In fact we have, 
\bgn{proposition}
If $\e(t)$ satisfies the equation $de^{\e(t)}\cdot\Ome=0$, then $\J_t$ is integrable.
\end{proposition}
\bgn{proof}
The equation (\ref{eq: de e(t) Ome=0}) is equivalent to the equation 
\bgn{equation}\label{eq: e -e(t) de e(t) Ome}
e^{-\e(t)}\, de^{\e(t)}\cdot\Ome =0
\end{equation}
Let $\pi_{U^{-n+3}}$ be the projection to the component $U^{-n+3}$. 
Then from \cite{Go1}, we have 
\bgn{align}\label{ali: -> M-C}
\pi_{U^{-n+3}}\( e^{-\e(t)}\, de^{\e(t)}\)\cdot\Ome =&\pi_{U^{-n+3}}d\Ome +
\pi_{U^{-n+3}}[d,\e(t)]\cdot\Ome+\frac12[\e(t), \e(t)]_S\Ome\\
=&\(\ol\pa\e(t)+\frac12[\e(t), \e(t)]_S\)\cdot\Ome =0
\end{align}
Then it follows that $\ol\pa\e(t)+\frac12[\e(t), \e(t)]_S=0$ on the complement $X\backslash D$. 
Since $\e(t)$ is a smooth section on $X$, we have $\ol\pa\e(t)+\frac12[\e(t), \e(t)]_S=0$
on $X$.
\end{proof}
\bgn{theorem}\label{th: unobstructedness  theorem}
Let $M=(X,J)$ be a compact K\"ahler manifold of dimension $n$.
We denote by $\J$ the \complex structure given by $J$.
We suppose that $M$ has an effective, anti-canonical divisor $D$. Then 
for every smooth global section 
$\e_1$ of $\w^2\ol \L(-D)\oplus \w^0 \ol \L(-D)$ with $d\e_1\cdot\Ome=0$,
there is a family of smooth global sections $\e(t)$ of $ \w^2\ol \L(-D)\oplus \w^0\ol \L(-D)$
such that 
$\J_t$ defined by $\e(t)$ is an integrable \complex structure and 
$\frac{d}{dt}\e(t)|_{t=0}=\e_1$ , where $t$ is a parameter of deformations which is sufficiently small.
\end{theorem}

{\indent\sc Proof of theorem \ref{th: unobstructedness  theorem}.}
We shall construct a family of smooth global sections 
${\displaystyle \e(t)=\e_1t+\e_2\frac{t^2}{2!}+\cdots}$ of  $\w^2\ol\L(-D)\oplus\w^0\ol \L(-D)$ for every  section 
$\e_1$ of $\w^2\ol \L(-D)\oplus\w^0\ol \L(-D)$ with $d\e_1\cdot\Ome=0$ which satisfies the following, 
\bgn{equation}\label{eq: 1-1}
de^{\e(t)}\cdot\Ome =0
\end{equation}
We denote by $\(e^{\e(t)}\)_{[i]}$ the $i$-th term of $e^{\e(t)}$ in $t$. 
Since both $\e_1\cdot\Ome $ and $\Ome$ are $d$-closed, we have 
$\(d e^{\e(t)} \)_{[1]}\cdot\Ome =\( e^{-\e(t)}de^{\e(t)}\)_{[1]}=d \e_1\cdot\Ome =0$.
We shall construct $\e(t)$ by the induction on $t$. 
We 
assume that there already exists a set of sections $\e_1,\cdots, \e_{k-1}$ of $\w^2\ol\L(-D)\oplus\w^0\ol \L(-D)$ such that 
\bgn{equation}\label{eq:1-4}
\( de^{\e(t)}\cdot\Ome \)_{[i]}=0, \qquad 0\leq \text{\rm for all } i<k
\end{equation}
The assumption (\ref{eq:1-4}) is equivalent to the following, 
\bgn{equation}\label{eq: 1-5}
\(e^{-\e(t)}\, d\, e^{\e(t)}\)_{[i]}\cdot\Ome=0,\qquad 0\leq\text{\rm for all } i<k
\end{equation}
Then $k$-th term is given by 
\bgn{equation}\label{eq: 1-6}
\(e^{-\e(t)}\, d\, e^{\e(t)}\)_{[k]}\cdot\Ome = \sum_{i+j =k}\( e^{-\e(t)}\)_{[i]}\, 
\( d e^{\e(t)}\)_{[j]}\cdot\Ome =\( de^{\e(t)}\)_{[k]}\cdot\Ome
\end{equation}
It follows from (\ref{ali: -> M-C}) that 
\bgn{equation}
\( e^{-e(t)}\, d\, e^{\e(t)}\)_{[k]}\cdot\Ome =\(d\e(t)\)_{[k]}\cdot\Ome +\frac12\([\e(t), \e(t)]_S\)_{[k]}\cdot\Ome
\end{equation}
We denote by $\Ob_k$ the non-linear term $\frac12\([\e(t), \e(t)]_S\)_{[k]}\cdot\Ome$. 
Since $[\e(t), \e(t)]_S$ is a section of $\w^3\ol L(-D)\oplus\w^1\ol L(-D)$, 
$\Ob_k$ is a section of $\bE^2=U^{-n+3}\oplus U^{-n+1}$. 
It follows from (\ref{eq: 1-6}) that $\Ob_k$ is a $d$-exact differential form. 
Hence $\Ob_k$ defines the cohomology class $[\Ob_k]\in H^2(\bE^\bullet)$ of the complex $(\bE^\bullet, d)$.
Since $M$ is a K\"ahler manifold, we apply the $\pa\ol\pa$-lemma to obtain the injective map $p^2$ from
$H^2(\bE^\bullet) $ to the direct sum of the de Rham cohomology groups. 
Since $\Ob_k$ is $d$-exact, the image of the class $p_2([\Ob_k])=0$. 
Hence the class $[\Ob_k]\in H^2(\bE^\bullet)$ vanishes, since the the map $p_2$ is injective.
Then the Hodge decomposition of the complex $(\bE^\bullet, d)$ shows that 
$\Ob_k=dd^*G(\Ob_k)$, where $d^*$ is the formal adjoint and $G$ is the Green operator of the complex $(\bE^\bullet, d)$.
Thus there is a unique section $\e_k$ of $\w^2\ol L(-D)\oplus \w^0\ol L(-D)$ such that 
$$
\frac1{k!}\e_k\cdot\Ome =-d^* G(\Ob_k)\in \bE^1,
$$
since $\w^2\ol L(-D)\oplus \w^0\ol L(-D)\cong \bE^1=U^{-n}\oplus U^{-n+2}$.
It follows that $\frac1{k!}d\e_k\cdot\Ome =-dd^* G(\Ob_k)=-\Ob_k.$
Then $\e_k$ satisfies the equation
$\( e^{-e(t)}\, d\, e^{\e(t)}\)_{[k]}=0$.
Thus by the induction, we obtain the power series $\e(t)$ which satisfies the equation (\ref{eq: 1-1}). 
As in \cite{Go1} the power 
 $\e(t)$ is a convergent series which is smooth.
\qed 
We shall obtain unobstructed K-deformations,  
\bgn{theorem}\label{th: unobstructedness  theorem-1} 
Let $M=(X,J)$ be a compact K\"ahler manifold of dimension $n$
and we denote by $\J$ the \complex structure given by $J$.
If $M$ has an effective, anti-canonical divisor $D$, 
then $M$ admits unobstructed $K$-deformations of \complex structures 
$\{ \J_t\}$ starting with $\J_0=\J$ which are 
parametrized by an open set of $H^2(\ol L(-D))\cong H^{n,2}\oplus H^{n-1,1}\oplus H^{n-2, 0}$ with the origin, that is, 
there is a family of smooth global sections $\e(t)$ of the sheaf $\w^2\ol\L(-D)$
such that 
$\J_t$ defined in (\ref{eq: Lt}) is an integrable \complex structure and 
$\frac{d}{dt}\e(t)|_{t=0}=\e_1$ for every representative 
$\e_1$ of $H^2(\ol L(-D))$ for small $t$, where $t$ is a parameter of deformations.
\end{theorem}
\bgn{proof}
Let $\e_1$ be a representative of the cohomology group $H^2(\ol L(-D))$.
Then $\e_1\cdot\Ome$ is a smooth differential form with $\ol\pa\e_1\cdot\Ome=0$.
Then it follows from the $\pa\ol\pa$-lemma that there is a function $k_1$ of $\w^0\L(-D)$ which satisfies 
$d(\e_1+\k_1)\cdot\Ome=0$. 
We put $\til\e_1=\e_1+\k_1$. 
Applying the theorem \ref{th: unobstructedness  theorem}, we obtain a section $\til\e(t)$ of 
$\w^2\ol\L(-D)\oplus\w^0\ol\L(-D)$ with 
$de^{\til\e(t)}\cdot\Ome=0$. 
The section $\til\e(t)$ is written as $\til\e(t) =\e(t)+\k(t)$, where 
$\e(t)\in \w^2\ol L(-D)$ and $\k\in \w^0\ol L(-D)$. Since $\Ad_{e^{\til\e(t)}}=\Ad_{e^{\e(t)}}$, 
the section $\displaystyle{\e(t)=\e_1t+\frac 1{2!}\e_2t^2+\cdots}$ gives $K$-deformations as we want.
\end{proof}

By taking $\e(t)$ as a family of global sections of $\w^2\ol L(-D)\cap (T^{1,0}\otimes\w^{0,1})$,
we have unobstructed deformations of usual complex structures $J_t$ which is given by the adjoint action $\Ad_{e^{\e(t)}}$, 
where $\Ad_{e^{\e(t)}}$ is a family of sections of GL$(TX, \C)$. Thus we obtain the following corollary, which is already obtained 
by Miyajima in the case of deformations of a normal isolated singularity \cite{Miya}.
\bgn{corollary}\label{cor: K-deformations of complex structures}
There is a family of deformations of complex structures $\{J_{t}\}$ starting with $J_{0}=J$ which satisfies 
$$
\frac{d}{dt}{\e(t)}|_{t=0} =\e_1,
$$
for every representative $\e_1$ of $H^1(M,\Theta(-D))$.
\end{corollary}
\bgn{proof}
We define a sheaf $\w^r\ol\L(-D)_{\SL}$ by 
$$
\w^r\ol\L(-D)_{\SL}(U)=\{\, f a\, |\, f\in I_D(U), \, a\in T^{1,0}\otimes\w^{0,r-1}\, \}
$$
The sheaf $\w^r\ol\L(-D)_{\SL}$ is the intersection $\w^r\ol \L(-D)\cap (T^{1,0}\otimes \w^{0,r-1})$ which is locally free. 
Then $\w^r\ol\L(-D)_{\SL}$ is a sheaf of smooth sections of the vector bundle $\w^r\ol L(-D)_{\SL}$.
As before, by the action on the meromorphic form $\Ome$ with a simple pole along the anti-canonical divisor $D$, we have the identification, 
$$
\w^r\ol L(-D)_{\SL}\cong \w^{n-1, r-1}
$$
Then we have the complex $(\w^\bullet\ol\L(-D)_{\SL}, \ol\pa)$ which is isomorphic to the Dolbeault complex $(\w^{n-1,\bullet}, \ol\pa)$.
We define vector bundles $\bE_{\SL}^\bullet$ by 
\bgn{align*}
&\bE^0_{\SL}=\w^{n-1,0}, \qquad\qquad \bE^1_{\SL}=\w^{n,0}\oplus \w^{n-1,1}\\
&\bE^2_{\SL}=\w^{n,1}\oplus\w^{n-1,2},  \quad\bE^3_{\SL}=\w^{n,2}\oplus\w^{n-1,3}, \cdots
\end{align*}
Then we have the complex $(\bE_{\SL}^\bullet, d)$ with the cohomology group $H^\bullet(\bE_{\SL})$,
$$
0\arrow \bE^0_{\SL}\overset{d}\arrow \bE^1_{\SL} \overset{d}\arrow\bE^2_{\SL}\overset{d}\arrow\cdots
$$
Let $\e_1$ be a representative of $H^1(M, \Theta(-D))$. Then $\e_1$ is a smooth global section of $\w^2\ol\L(-D)_{\SL}$ with $\ol\pa\e_1\cdot\Ome=0$. 
Then it follows from the $\pa\ol\pa$-lemma that there is a section $\k_1$ of  $\w^0\ol\L(-D)$ with 
$d(\e_1+\k_1)\cdot\Ome=0$. We put $\til\e_1=\e_1+\k_1$.
For such a section $\til\e_1$ of $\w^2\ol\L(-D)_{\SL}\oplus \w^0\ol\L(-D)$ with $d\til\e_1\cdot\Ome =0$, we shall construct a family of sections $\til\e(t)$ of $\w^2\ol\L(-D)_{\SL}\oplus \w^0\ol\L(-D)$ which satisfies 
$$
de^{\til\e(t)}\cdot\Ome=0,
$$
where $\til\e(t)=\til\e_1t+\frac1{2!}\til\e_2t^2+\cdots$.
As in the proof of the theorem \ref{th: unobstructedness  theorem}, we have the obstruction $(\Ob_k)_{\SL}$ which gives the class $[(\Ob_k)_{\SL}]\in H^2(\bE_{\SL})$. 
It suffices to show that the class $[(\Ob_k)_{\SL}]$ vanishes. 
In fact, the $(\Ob_k)_{\SL}\in \bE^2_{\SL}$ is a $d$-exact differential form. 
It follows from the Hodge decomposition that the map $p^\bullet_{\SL}$ from $H^\bullet(\bE_{\SL})$ to the direct sum of the de Rham cohomology groups is injective. 
Since $(\Ob_k)_{\SL}$ is $d$-exact, the image $p^2_{\SL}([(\Ob_k)_{\SL}])$ vanishes. 
Thus the class $[(\Ob_k)_{\SL}]\in H^2(\bE_{\SL})$ vanishes also.
Hence we have $\til\e(t)$ with $de^{\til\e(t)}\cdot\Ome=0$. 
Let $\e(t)$ be the component of $\til\e(t)$ of $\w^2\ol\L(-D)_{\SL}$. 
Since $\Ad_{e^{\til\e(t)}}=\Ad_{e^{\e(t)}}$, 
we have deformations of \complex structures $\J_t$ given by $\e(t)$.
\end{proof}
In the cases of K\"ahler surfaces, we obtain the following,
\bgn{corollary}\label{cor: unobstructedness S}
Let $S$ be a compact K\"ahler surface with the complex structure $J$ and 
a K\"ahler form $\ome$. 
If $S$ has an effective, anti-canonical divisor $[D]=-K_S$, 
then $S$ admits unobstructed deformations of \complex structures parametrized by an open set of 
the full cohomology group $H^0(S)\oplus H^2(S)\oplus H^4(S)$ of even degree
on $S$ 
\end{corollary}
{\indent\sc Proof of Corollary \ref{cor: unobstructedness S}.} 
Since $-K_S$ is effective, then we have the vanishing $H^2(S, {\cal O}_S) =
H^0(S, K_S)\cong H^0(S, I_D)=\{0\}$.
Thus $H^2(S)\cong H^{1,1}$. Then the result follows from 
the theorem \ref{th: unobstructedness  theorem}.
\qed
 
 \section{K-deformations of \complex structures in terms of $\CL^2(-D)$}
 We denote by $\w^r L(-\ol D)$ the complex conjugate of the bundle $\w^r\ol L(-D)$. 
 A section of $\w^r L(-\ol D)$ is locally written as 
 $\ol f\ol a$ for $f\in I_D$ and $a\in \w^r\ol\L$.
 Let $\(\w^2\ol L(-D)\oplus \w^2L(-\ol D)\)^\R$ be the real part of the bundle
 $\w^2\ol L(-D)\oplus \w^2L(-\ol D)$, which is the subbundle of $\CL^2$.
We define a bundle $\CL^2(-D)$ by 
$$\CL^2(-D):=\(\w^2\ol L(-D)\oplus \w^2L(-\ol D)\)^\R\oplus  \w^0 \ol L(-D)$$
 \bgn{lemma}\label{lem: e(t)-> a(t)}
For small deformations of almost \complex structures
$\J_t$  given by  a family of smooth global sections $\e(t)$ of $ \w^2\ol L(-D)\oplus \w^0\ol L(-D)$ as in (\ref{eq: Lt}), there exists a unique family of global sections $a(t)$ of the bundle CL$^2(-D)$ such that
$$
e^{\e(t)}\cdot \Ome =e^{a(t)}\cdot\Ome.
$$
that is,
$\J_t=\Ad_{e^{a(t)}}\J_0$. 

Conversely if we have a family of deformations of almost \complex structure $\J_t
=\Ad_{e^{a(t)}}\J_0$ 
which is given by the action of a family of global sections $a(t)$ of $\CL^2(-D)$, 
then there exists a unique family of global sections $\e(t)$ of $\w^2\ol L(-D)\oplus\w^0 \ol L(-D)$ such that 
$\J_t$ is given by the action of $\e(t)$ and $\e^{a(t)}\cdot\Ome =
e^{\e(t)}\cdot\Ome$.
\end{lemma}
\bgn{proof}
For a section $\e$ of  $\w^2\ol L(-D)\oplus \w^0\ol L(-D)$, 
we have a unique $a\in \Gam(X, \CL^2(-D))$ such that 
$e^{\e}\cdot\Ome =e^a\cdot\Ome$. 
Conversely, there is a unique section $\e$ of $ \w^2\ol L(-D)\oplus \w^0\ol L(-D)$ such that
$e^{\e}\cdot\Ome =e^a\cdot\Ome$ for any section $a$ of $\CL^2(-D)$. 
Then applying the method in \cite{Go1}, we obtain the result.
\end{proof}
The operator $e^{-a(t)}\circ d\circ e^{a(t)}$ acting on $K_J=U^{-n}$ is 
already discussed in \cite{Go0} which  
is a Clifford-Lie operator of order $3$ whose image is in $U^{-n+1}\oplus U^{-n+3}$.

It is shown in \cite{Go1} that the almost \complex structure $\J_t=\Ad_{e^{a(t)}}\J$ is integrable if and only if 
the projection to the component $U^{-n+3}$ vanishes, that is,
$$
\pi_{U^{-n+3}}e^{-a(t)}\circ d\circ e^{a(t)}\cdot\Ome=0
$$
In particular, $e^{-a(t)}\circ d\circ e^{a(t)}\cdot\Ome=0$ implies that the $\J_t$ is integrable.
We denote by $\( e^{-a(t)}\circ d\circ e^{a(t)}\)_{[k]}$ the $k$ th term of $e^{-a(t)}\circ d\circ e^{a(t)}$.

Thus by the theorem \ref{th: unobstructedness  theorem} and the lemma 
\ref{lem: e(t)-> a(t)}, we have the following,

\bgn{proposition}\label{prop: k th-G.complex}
Let $M=(X, J)$ be a compact K\"ahler manifold with a K\"ahler form $\ome$.
We assume that $M$ has an effective, anti-canonical divisor $D$.
If there is a set of global sections $a_1,\cdots, a_{k-1}$ of CL$^2(-D)$ which satisfies 
\bgn{equation}\label{eq: k th-G.complex}
\(e^{-a(t)}\, de^{a(t)}\)_{[i]}\cdot\Ome=0, \qquad 0\leq\text{\rm for all }
i<k,
\end{equation}
and $\|a(t)\|_s\underset{k-1} \<C_1M(t)$,
then there is a global section $a_k$ of CL$^2(-D)$ which satisfies the followings:
$$
\(e^{-a(t)}\, de^{a(t)}\)_{[k]}\cdot\Ome=0
$$
and $\|a(t)\|_s\underset{k}\< C_1\lam M(t)$, where 
$a(t) =\sum_{i=1}^\infty \frac 1{i!} a_i t^i$ and 
$M(t)$ is the convergent series (\ref{eq: M(t)}) in section 7  and
$C_1$ is a positive constant and 
$\|a(t)\|_s$ denotes the Sobolev norm of $a(t)$. 
\end{proposition}
 (The proof of the inequality $\|a(t)\|_s\underset{k}\< C_1\lam M(t)$ is already seen in \cite{Go1}, see proposition 1.1 and 1.4  in \cite{Go3} for more detail).

\section{Deformations of generalized K\"ahler structures}
Let $(X, J,\ome)$ be a compact K\"ahler manifold with an effective anti-canonical divisor $D$
and $(\J, \J_{\psi})$ the generalized K\"ahler structure induced from
$(J,\ome)$ by $\J=\J_J$ and $\psi=e^{\sqrt{-1}\ome}$.
Since two \complex structures $\J$ and $ \J_{\psi}$ are commutative, 
the generalized K\"ahler structure $(\J, \J_{\psi})$ gives the simultaneous  decomposition of 
$(\TT)^\C$, 
$$
(\TT)^\C=L^+_\J\oplus L^-_\J\oplus \ol L^+_\J\oplus\ol L^-_\J,
$$
where $L^+_\J\oplus L^-_\J$ is the eigenspace with eigenvalue $\sqrt{-1}$ with respect to $\J$ and $L^+_\J\oplus \ol L^-_\J$ is the eigenspace with eigenvalue $\sqrt{-1}$ with respect to $\J_\psi$ and $\ol L^\pm_\J$ denotes the complex conjugate.
In \cite{Go1,Go2}, the author showed the stability theorem of 
generalized K\"ahler structures with one pure spinor, which implies that 
if there is a one dimensional  analytic deformations of \complex structures $\{\J_t\}$ parametrized by $t$,
then there exists a family of non-degenerate, $d$-closed pure spinor $\psi_t$
such that the family of pairs $(\J_t, \psi_t)$ becomes deformations of 
{generalized K\"ahler} structures starting from $(\J, \psi)=(\J_0, \psi_0)$.
As in section 2, small K-deformations $\J_t$ are given by the adjoint action of a family of sections $a(t)$ of CL$^2(-D)$, 
$$
\J_t:=\Ad_{e^{a(t)}}\J_0. 
$$
Then we can obtain a family of real sections $b(t)$ of the bundle
$(L^-_{\J_0}\cdot\ol L^+_{\J_0}\oplus \ol L^-_{\J_0}\cdot L^+_{\J_0})^\R$ such that $\psi_t=e^{a(t)}e^{b(t)}\psi_0$ is a family of non-degenerate, $d$-closed pure spinor.  
The bundle $K^1=U^{0,-n+2}$  is generated by the action of real sections of $(L^-_{\J_0}\cdot\ol L^+_{\J_0}\oplus \ol L^-_{\J_0}\cdot L^+_{\J_0})$ on $\psi$ (see page 125 in \cite{Go2} for more detail).

We define a family of sections $Z(t)$ of $\CL^2$ by 
$$
e^{Z(t)}= e^{a(t)}\, e^{b(t)}.
$$
Since $\Ad_{e^b(t)}\J_0=\J_0$, we obtain 
$\J_t=\Ad_{e^{a(t)}}\J_0=\Ad_{e^{a(t)}}\, \Ad_{e^{b(t)}}\J_0=\Ad_{e^{Z(t)}}\J_0 $. Then the family of deformations of generalized K\"ahler structures is given by the action of $e^{Z(t)}$, 
$$
(\J_t, \,\psi_t) =\(\Ad_{e^{Z(t)}}\J_0, \, \,\,e^{Z(t)}\cdot\psi \,\).
$$
By the similar method as in \cite{Go2} together with the proposition \ref{prop: k th-G.complex}, we obtain the following proposition, 
\bgn{proposition}\label{prop: GK deformations}
Let $M=(X, J)$ be a compact K\"ahler manifold with a K\"ahler form $\ome$.
We assume that $M=(X, J)$ has an anti-canonical divisor $D$.
If there is a set of sections $a_1,\cdots a_{k-1}$ of CL$^2(-D)$ which satisfies 
$$
\(e^{-a(t)}\, de^{a(t)}\)_{[i]}\cdot\Ome=0, \qquad 0\leq \text{\rm for all }
i<k,
$$
and $\|a(t)\|_s \underset{k-1}\< K_1M(t)$ for a positive constant $K_1$,
then there is a set of real sections $b_1,\cdots, 
b_k$ of the bundle 
$(L^-_{\J_0}\cdot\ol L^+_{\J_0}\oplus \ol L^-_{\J_0}\cdot L^+_{\J_0})$ which satisfies 
the following equations:
\bgn{align}
&\(e^{-Z(t)}\, de^{Z(t)}\)_{[k]}\cdot\Ome=0\\
&\( d e^{Z(t)}\cdot\psi_0\)_{[i]}=0,\qquad \text{\rm for all } i\leq k\\
&\|a(t)\|_s\underset {k}\< K_1\lam M(t)\\
&\|b(t)\|_s\underset{k}\<K_2M(t)
\end{align}
where $a_k$ is the section constructed in the proposition \ref{prop: k th-G.complex} and
$e^{Z(t)}=e^{a(t)}\, e^{b(t)}$ and $M(t)$ is the convergent series in the proposition \ref{prop: k th-G.complex} and 
a positive constant $K_2$ is determined by $\lam$ and $K_1$.
The constant $\lam$ in $M(t)$ is  sufficiently small which will be suitably selected to show the convergence of the power series $Z(t)$ as in \cite{Go1}.
\end{proposition}
\section{Deformations of bihermitian structures}
There is a one to one correspondence between generalized K\"ahler structures and bihermitian structures with 
 the condition (\ref{eq: torsion condition}) \cite{Gu1}.
 In this section we shall give an explicit description of $\Gam^\pm_t$ which gives rise to deformations of 
 bihermitian structures $(J^+_t, J^-_t)$ 
 corresponding to deformations of generalized K\"ahler structures with one pure spinor $(\J_t, \psi_t)$ in section 4.
 The correspondence is defined at each point on a manifold. 
 The non-degenerate, pure spinor $\psi_t$ induces the \complex structure $\J_{\psi_t}$. 
 Since $(\J_t, \J_{\psi_t})$ is a generalized K\"ahler structure and $\J_t$ commutes with $\J_{\psi_t}$, 
 we have the simultaneous decomposition of $(\TT)^\C$ into four eigenspaces as before, 
 $$
 (\TT)^\C=L^+_{\J_t}\oplus L^-_{\J_t}\oplus \ol {L_{\J_t}^+}\oplus \ol {L_{\J_t}^-},
 $$
where each eigenspace is given by the intersection of eigenspaces of both $\J_t$ and $\J_{\psi_t}$, 
 \bgn{align*}
 &L_{\J_t}^-=L_{\J_t}\cap \ol L_{\psi_t},\qquad \ol {L_{\J_t}^+}=\ol L_{\J_t}\cap \ol L_{\psi_t}\\
 &L_{\J_t}^+=L_{\J_t}\cap L_{\psi_t},\qquad\ol {L_{\J_t}^-}=\ol L_{\J_t}\cap L_{\psi_t},
   \end{align*}
where $ L_{\J_t}$ is the eigenspace of $\J_t$ with eigenvalue $\sqrt{-1}$
and $L_{\psi_t}$ denotes the eigenspace of $\J_{\psi_t}$ with eigenvalue $\sqrt{-1}$.
Since $\J_t=\Ad_{e^{Z(t)}}(\J_0)=\Ad_{e^{Z(t)}}\circ \J_0\circ\Ad_{e^{-Z(t)}}$ and 
$\J_{\psi_t}=\Ad_{e^{Z(t)}}(\J_\ome)$, we have the isomorphism between eigenspaces,
  $$
  \Ad_{e^{Z(t)}} : \ol{L_{\J_0}^\pm}\to \ol{L_{\J_t}^\pm}.
  $$
Let $\pi$ be the projection from $\TT$ to the tangent bundle $T$. We restrict the map $\pi$ to the 
eigenspace $\ol{L_{\J_t}^\pm}$ which yields the map $\pi_t^\pm : \ol{L_{\J_t}^\pm}\to T^\C$.
Let $T^{1,0}_{J^\pm_t}$ be the complex tangent space of type $(1,0)$ with respect to 
  $J^\pm_t$. 
  Then it follows that $T^{1,0}_{J^\pm_t}$ is given by the image of $\pi_t^\pm$, 
  $$
  T^{1,0}_{J^\pm_t}=\pi_t^\pm( \ol{L^\pm_{\J_t}})
  $$
Since deformations of generalized K\"ahler structures are given by the action of $e^{Z(t)}$, 
the ones of bihermitian structures $J^\pm_t$ should be described by the action of $\Gam^\pm_t$ of the bundle GL$(T)$ 
which is obtained from $Z(t)$. 
We shall describe $\Gam^\pm_t$ in terms of $a(t)$ and $b(t)$.
A local basis of $\ol{L_{\J_0}^\pm}$ is given by 
$$\{ \Ad_{e^{\pm\sqrt{-1}\ome}}V_i= V_i\pm \sqrt{-1}[\ome, V_i]\, \}_{i=1}^n,$$
for a local basis $\{ V_i\}_{i=1}^n$ of $T^{1,0}_J$, 
where we regard $\ome$ as an element of the Clifford algebra and then the bracket 
$[\ome, V_i]$ coincides with the interior product $i_{V_i}\ome$.
It follows that the inverse map $(\pi^\pm_0)^{-1} : T^{1,0}_{J} \to \ol {L_{\J_0}^\pm}$ 
is given by the adjoint action of $e^{\pm\sqrt{-1}\ome}$, 
\bgn{equation}\label{eq: pi0{-1}}
\Ad_{e^{\pm\sqrt{-1}\ome}}=(\pi^\pm_0)^{-1}.
\end{equation}
 We define a map $(\Gam_t^\pm)^{1,0}: T^{1,0}_J \to T^{1,0}_{J^\pm_t}$ by the composition, 
 \bgn{align}
 (\Gam_t^\pm)^{1,0}=&\pi^\pm_t\circ \Ad_{e^{Z(t)}}\circ (\pi^\pm_0)^{-1}\\
 =&\pi\circ \Ad_{e^{Z(t)}}\circ \Ad_{e^{\pm\sqrt{-1}\ome}}
\end{align}
\medskip
$$\xymatrix{
\ol{L^\pm_{{\cal J}_0}}\ar@{->}[rr]^{\text{\rm Ad}_{e^{Z(t)}} }\ar@{->}[d]_{\pi^\pm_0}& &\ol{L^\pm_{{\cal J}_t}}\ar@{->}[d]^{\pi^\pm_t}\\
T_J^{1,0}\ar@{->}[rr]_{(\Gamma^\pm_t)^{1,0}}&& T^{1,0}_{J_t^\pm}
}
$$

Together with the complex conjugate $(\Gam^\pm_t)^{0,1} :T^{0,1}_J \to T^{0,1}_{J^\pm_t}$, 
we obtain the map $\Gam^\pm_t$ which satisfies 
$J^\pm_t= (\Gam_t^\pm)^{-1} \circ J\circ \Gam_t^\pm$. 

Let $J^*$ be the complex structure on the cotangent space $T^*$ which is given by 
$\lan J^*\eta,  v\ran=\lan \eta, Jv\ran$, where $\eta\in T^*$ and $v\in T$ and 
$\lan \,,\,\ran$ denote the coupling between $T$ and $T^*$.
We define a map $\h J^\pm : \TT \to \TT$ by $\h J^\pm
(v,\eta) = v\mp J^*\eta$ for $v\in T$ and $\eta\in T^*$.
Then $\Gam^\pm_t$ is written as 
\bgn{align}\label{eq: Gam t pm }
\Gam^\pm_t= &\pi\circ \Ad_{e^{Z(t)}}\circ \h J^\pm\circ\Ad_{e^\ome}\\
=&\pi\circ\Ad_{e^{a(t)}}\circ \Ad_{e^b(t)}\circ \h J^\pm\circ\Ad_{e^\ome}\in\text{\rm End}(T),
\end{align}
where note that $\h J^\pm\circ\Ad_{e^{\ome}}(T^{1,0}_J)=\ol{L^\pm_{{\cal J}_0}}.$
The $k$ th term of $\Gam^\pm_t$ is denoted by $(\Gam_t^\pm)_{[k]}$ as before. 
Note that $(\Gam_t^\pm)_{[0]}=$id$_T$.
We also put $\Gam^\pm (a(t), b(t))=\Gam_t^\pm$.
\bgn{lemma}\label{lem: Gam-k}
The $k$ th term $(\Gam_t^\pm)_{[k]}$ is given by 
\bgn{align*}
(\Gam_t^\pm)_{[k]}= \frac 1{k!}\pi\circ(\ad_{a_k}+\ad_{b_k})\circ\h J^\pm\circ\Ad_{e^\ome}+\wtil{\Gam^\pm_k}(a_{<k}, b_{<k})
\end{align*}
where the second term $\wtil{\Gam^\pm_k}(a_{<k}, b_{<k})$ depends only on $a_1,\cdots, a_{k-1}$ and 
$b_1,\cdots, b_{k-1}$.
\end{lemma}
\bgn{proof}
Substituting the identity $\Ad_{e^{Z(t)}}=$id$ +\ad_{Z(t)}+\frac1{2!}(\ad_{Z(t)})^2+\cdots$, 
we have 
\bgn{align}
\Gam_t^\pm =&\pi\circ\Ad_{e^{Z(t)}}\circ \h J^\pm\circ\Ad_{e^\ome}\\
=&\pi\circ\(\sum_{i=0}^\infty \frac 1{i!}\ad^i_{Z(t)}\circ \h J^\pm\circ \Ad_{e^\ome}\)\\
\end{align}
Then $k$-th term is given by 
\bgn{align}
\(\Gam_t^\pm \)_{[k]}
=&\pi\circ\(\ad_{Z(t)}\circ\h J^\pm\circ\Ad_{e^\ome}\)_{[k]}+\sum_{i=2}^k\pi\circ\(\frac 1{i!}(\ad_{Z(t)}^i\circ \h J^\pm\circ\Ad_{e^\ome}\)_{[k]}\\
=&\frac 1{k!}\pi\circ(\ad_{a_k}+\ad_{b_k})\circ\h J^\pm\circ\Ad_{e^\ome}+\wtil{\Gam^\pm_k}(a_{<k}, b_{<k}),
\end{align}
where $\wtil{\Gam^\pm_k}(a_{<k}, b_{<k})$ denotes the non-linear term depending $a_1,\cdots, a_{k-1}$ and $b_1, \cdots, b_{k-1}$.
\end{proof}
\bgn{lemma}\label{lem: key b }
Let $b$ be a section of the bundle 
$(L^-_{\J}\cdot\ol L^+_{\J}\oplus \ol L^-_{\J}\cdot L^+_{\J})$. 
Then we have 
$$
[\,\pi\circ\Ad_{e^b}\circ\h J^\pm\circ\Ad_{e^\ome},\, J\, ]=0\in\text{\rm End}(T).
$$
\end{lemma}
\bgn{proof}
For simplicity, we write $L^\pm$ for $L_{\J}^\pm$.
Note that $[\pi\circ\h J^\pm\circ\Ad_{e^\ome}, \, J]=[\id_T, J]=0$ and 
$\Ad_{e^b}=\id+\ad_b+\frac 1{2!}(\ad_b)^2+\cdots.$
 We also recall the image $\h J^\pm\circ\Ad_{e^\ome}(T^{1,0}_J) =(\pi^{\pm}_0)^{-1}(T^{1,0}_J)=\ol L^\pm_J$.
Since $b$ is a section $\ol L^+\cdot L^-\oplus L^+\cdot\ol L^-$, 
the image $(\ad_b)^n(\ol L^\pm)$ is given by 
$$\bgn{cases}
(\ad_b)^n(\ol L^\pm)=& \ol L^\pm, \quad (n : even)\\
(\ad_b)^n(\ol L^\pm)=& \ol L^\mp, \quad (n : odd)\\
\end{cases}
$$
Since $\pi(\ol L^\pm)=T^{1,0}_J$, we have 
$\pi\circ(\ad_b)^n\circ\h J^\pm\circ\Ad_{e^\ome}(T^{1,0}_J)=T^{1,0}_J$. 
Thus $\pi\circ(\ad_b)^n\circ\h J^\pm\circ\Ad_{e^\ome}\in \End(T)$ preserves $T^{1,0}_J$. Hence we have
$$
[\pi\circ(\ad_b)^n\circ\h J^\pm\circ\Ad_{e^\ome}, \, J]=0.
$$
Then the  result follows.
\end{proof}

The tensor space $T\otimes T^*$ defines a subbundle of CL$^2$. We denote it by $T\cdot T^*$.
An element $\gam\in T\cdot T^*$ gives the endmorphism ad$_\gam$ by ad$_\gam E=[\gam, E]$ for $E\in \TT$, which preserves the tangent bundle
$T$ and the cotangent bundle $T^*$ respectively. 
We also regard ad$_\gam$ as a section of End$(T)$. 
\bgn{lemma}\label{lem: gam}
Let $\gam$ be an element of $T\cdot T^*$. Then we have 
$$\pi\circ(\ad_{\gam}\circ \h J^\pm\circ \Ad_{e^\ome})=\ad_\gam\in \text{\rm End}(T). $$
\end{lemma}
\bgn{proof}
For a tangent vector $v\in T$, we have $\Ad_{e^\ome}v= v+[\ome, v]=v+\ad_\ome v$.
Since the map $\ad_\gam$ preserves the cotangent $T^*$, 
we have $\ad_{\gam}\circ \h J^\pm\circ \ad_{\ome}(v)\in T^*$ for all tangent $v\in T$.
Thus it follows that
$\pi( \ad_\gam\circ\h J^\pm\circ\ad_\ome)=0$, since $\pi$ is the projection to the tangent $T$.
Thus we obtain the result.
\end{proof}
\bgn{lemma}\label{lem: d Gam-k }
We assume that there is a set of sections 
$a_1,\cdots , a_k$ of $\CL^2(-D)$ and real sections $b_1,\cdots, b_k$ of 
$(L^-_{\J_0}\cdot\ol L^+_{\J_0}\oplus \ol L^-_{\J_0}\cdot L^+_{\J_0})$ which satisfies the following equations, 
\bgn{align*}
&\( e^{-Z(t)}\, d\, e^{Z(t)}\)_{[i]}\cdot\Ome=0, \quad 0\leq 
\forall\, i\leq k\\
&\( d e^{Z(t)}\cdot \psi_0 \)_{[i]}=0, \,\,\,\,\,\quad 0\leq 
\forall\, i\leq k\\
&[(\Gam^\pm_t)_{[i]}, \, J]=0,  \,\,\,\,\,\,\,\, \quad\quad 0\leq 
\forall\, i< k\\
\end{align*}
Then the $k$-th term $(\Gam_t^\pm)_{[k]}$ satisfies 
$$
\pi_{U^{-n+3}}[\, d,\,(\Gam_t^\pm)_{[k]}\, ]  =0,
$$
where $[\, d,\,(\Gam_t^\pm)_{[k]}\, ] $ is an operator from $U^{-n}=K_\J$ to 
$U^{-n+1}\oplus U^{-n+3}$ and $\pi_{U^{-n+3}}$ denotes the projection to the component $U^{-n+3}$.
\end{lemma}

\bgn{proof} 
Since the obstructions to $K$-deformations of \complex structures vanishes, 
we obtain a family of section $\ch a(t)$ with $\ch a_i =a_i$ for $ i=1, \cdots k$ such that $\ch a(t)$ gives $K$-deformations of \complex structures, that is, 
$$
\pi_{U^{-n+3}}e^{-\ch a(t)} d e^{\ch a(t)} \cdot\Ome=0.
$$
The the stability theorem of generalized K\"ahler structures in \cite{Go1} provides deformations of generalized K\"ahler structures with one pure spinor:
$$(\Ad_{e^{\ch Z(t)}}\J_0, \,e^{\ch Z(t)}\psi_0),$$where $e^{\ch Z(t)}=e^{\ch a(t)}e^{\ch b(t)}$, where 
$\ch b(t)$ is a family of real sections with $\ch b_i =b_i$, for $i=1,\cdots, k$.
From the correspondence between generalized K\"ahler structures and bihermitian structures, we  have the family of bihermitian structures 
$(J^+_t, J^-_t)$ which is given by the action of $\ch{\Gam}^\pm_t:=\Gam^\pm_t(\ch a(t), \ch b(t))$ of GL$(T)$. 
Since $J_t^\pm$ is integrable, we have 
\bgn{equation}\label{eq: ch Gam}
\pi_{U^{-n+3}}\((\ch{\Gam}^\pm_t)^{-1} \, d\,\ch{\Gam}^\pm_t\)=0.
\end{equation}
Let $\Ome$ be a $d$-closed meromorphic form of type $(n,0)$ with a simple pole along $D$ as before.
Then we have 
$$
d\Gam^\pm_t\Ome \underset{k}\equiv \Gam^\pm_t E(t) \Ome.
$$
Since $d\Ome=0$, the degree of $E(t)$ is greater than or equal to $1$.
The condition $[(\Gam^\pm_t)_{[i]}, \, J]=0$ $( 0\leq i<k)$
implies that $(\Gam^\pm_t)_{[i]}E(t) \Ome  \in U^{-n+1}_\J$. 
Thus we have 
$$
d (\Gam^\pm_t)_{[k]}\Ome=\sum_{\stackrel {\scriptstyle i+j=k}{ 0<i,j<k}} (\Gam^\pm_t)_{[i]} E(t)_{[j]}\Ome\in U^{-n+1}_\J
$$
Hence we have 
$\pi_{U^{-n+3}}[d, (\Gam^\pm_t)_{[k]}]=0$.
\end{proof}

\bgn{lemma}\label{lem: ad gam}
For a section $a$ of $\CL^2(-D)$ and every section $P$ of $\End(\TT)$, we define a section
$\zeta$ of $T\cdot T^*$ by 
$$
\ad_\zeta =[\pi \circ\ad_a\circ P|_T, \, J]
$$
Then $\zeta$ is a section of $\CL^2(-D)$, where $P|_T : T\to \TT$ denotes the restriction to the tangent bundle $T$.
\end{lemma}
\bgn{proof}
As in section2, we have 
$$ \w^2\ol L(-D)=(\w^{0,2}\otimes[-D])
\oplus(T^{1,0}\otimes\w^{0,1}\otimes[-D])\otimes(T^{2,0}\otimes[-D]),$$
Thus a section $\e\in \w^2\ol L(-D)$ gives 
$\pi\circ\ad_\e(E)\in T^{1,0}(-D)$ for all $E \in \TT$.
Since CL$^2(-D)=(\w^2\ol L(-D)\oplus{\w^2 L(-\ol D)})^\R\oplus\w^0\ol L(-D)$, we have $\pi\circ\ad_a(E) \in 
T^{1,0}(-D)\oplus T^{0,1}(-\ol D)$ for all $E\in\TT$.
Hence $\pi\circ\ad_a\circ P|_T$ is a section of 
$\(T^{1,0}(-D)\oplus T^{0,1}(-\ol D)\)\otimes\w^1.$
Taking the bracket, it turns out that 
$[ \pi\circ\ad_a\circ P|_T, J]$ is a section of  $(T^{1,0}(-D)\otimes\w^{0,1})
\oplus (T^{0,1}(-\ol D)\otimes\w^{1,0})$ which is the subbundle $\w^2\ol L(-D)\oplus \w^2 L(-\ol D)$. 
Thus $\zeta$ is a real section of $\w^2\ol L(-D)\oplus \w^2 L(-\ol D)
\subset \CL^2(-D)$. 
\end{proof}
\bgn{lemma}\label{lem: ad gam k}
We define a section $\zeta_k$ of  $T\cdot T^*$ by 
$$
\ad_{\zeta_k} =\big[(\Gam^+_t)_{[k]} , \, J\big] \in \End(T)
$$
for a section $a(t)$ of $\CL^2(-D)$ and $b(t)\in(L^-_{\J}\cdot\ol L^+_{\J}\oplus \ol L^-_{\J}\cdot L^+_{\J})$, 
where $(\Gam^+_t)_{[k]}:=\Gam^+(a(t),b(t))_{[k]}.$
Then $\zeta_k$ is a section of $\CL^2(-D)$.
Further we define $\gam_k$ by 
$$
\ad_{\gam_k}:=\frac{-1}{(2\sqrt{-1})^2}[\ad_{\zeta_k}, J].
$$
Then $\gam_k\in T\cdot T^*$ is also a section of $\CL^2(-D)$ which satisfies 
$$
\ad_{\zeta_k}+[\ad_{\gam_k}, J]=0.
$$
\end{lemma}
\bgn{proof}
From the description in (\ref{eq: Gam t pm }), 
\bgn{align}
\Gam^+_t
=&\pi\circ\Ad_{e^{a(t)}}\circ \Ad_{e^{b(t)}}\circ \h J^+\circ\Ad_{e^\ome}\\
=&\pi\circ(\Ad_{e^{a(t)}}-\id)\circ Q+ \pi\circ Q
\end{align}
where $Q:= \Ad_{e^{b(t)}}\circ \h J^+\circ\Ad_{e^\ome}$.

From the lemma \ref{lem: key b }, we have 
$[\pi\circ Q, J]=0$. Thus we have 
$$
[\Gam^+_t, J]= [\pi\circ(\Ad_{e^{a(t)}}-\id)\circ Q|_T,\,\, J]
$$
Since we have  $\Ad_{e^{a(t)}}-\id=\ad_{a(t)}\circ R$, 
where $\displaystyle{R= \sum_{j=1}^\infty\frac1{j!}\ad_{a(t)}^{j-1}}$.
If we set $P=R\circ Q$, we have 
$$
[\Gam^+_t, J] = [\pi\circ \ad_{a(t)}\circ P|_T,\,\, J]
$$
Then it follows from the lemma \ref{lem: ad gam} that $\zeta_k$ is a section of $\CL^2(-D)$. 
We decompose $\ad_{\zeta_k}$ by 
$$\ad_{\zeta_k}=(\ad_{\zeta_k})'+(\ad_{\zeta_k})'', $$ where $(\ad_{\zeta_k})'\in T^{1,0}(-D)\otimes\w^{0,1}$ and 
$(\ad_{\zeta_k})''\in T^{0,1}(-\ol D)\otimes \w^{1,0}$.
Then the bracket is given by 
$[\ad_{\zeta_k}, J]= -2\sqrt{-1}(\ad_{\zeta_k})'+2\sqrt{-1}(\ad_{\zeta_k})''\in \CL^2(-D)$.
Thus $\gam_k$ is also a section of $\CL^2(-D)$ which satisfies 
$$
\ad_{\zeta_k}+[\ad_{\gam_k}, J]=0.
$$
\end{proof}
\bgn{lemma}\label{lem: rho k}
Let $\Gam^+_t$ be a section of GL$(T)$ given in the lemma \ref{lem: d Gam-k } and 
$\zeta_k$ and $\gam_k$ be as in the lemma \ref{lem: ad gam k}.  
Then there is a global function $\rho_k$ of $\w^0\ol \L(-D)$ such that 
$$
d\gam_k\cdot\Ome =d(\rho_k\Ome),
$$
\end{lemma}
\bgn{proof}
The condition $\pi_{U^{-n+3}}[\, d,\,(\Gam_t^+)_{[k]}\,]  =0$ in the lemma \ref{lem: d Gam-k } implies that 
$(d\Gam_t^\pm)_{[k]}\cdot\Ome \in U^{-n+1}$.
Thus we have that $d\gam_k\cdot\Ome\in U^{-n+1}$.
Since $\ad_{\gam_k}$ is a section of GL$(T)$, we see that
$d\gam_k\cdot\Ome$ is a $d$-exact form of type $(n,1)$. 
Then applying the $\pa\ol\pa$-lemma, it turns out that 
$ d\gam_k\cdot\Ome=d\rho_k\Ome$ for a smooth function $\rho_k$.
Since $\rho_k\Ome$ is smooth, we have $\rho_k$ is a global function of $\w^0\ol\L(-D)$.
\end{proof}


\section{Bihermitian structures on compact K\"ahler surfaces}
Let $S$ be a compact K\"ahler surface with a K\"ahler form $\ome$ with 
an anti-canonical divisor $D$. 
The divisor $D$ is given as the zero locus of a section $\b\in H^0(S, K^{-1})\cong H^0(S, \w^2\Theta)$. Then the section $\b$ is also regarded as a section of 
$H^0(S, \w^2\Theta(-D))$ which is a holomorphic Poisson structure vanishing along the divisor $D$.
The contraction $\b\cdot\ome$ of $\b$ by $\ome$ is defined by the commutator $[\b, \ome]$ which is a $\ol\pa$-closed $T^{1,0}$-valued form of type $(0,1)$. 
Let $\Ome$ be the meromorphic $2$-form on $S$ with a pole along the divisor $D$ with $\b \cdot \Ome =1$. 
Then we have 
$$(\b\cdot\ome)\cdot\Ome =[\b,  \ome]\cdot\Ome =-\ome,
$$
since $\ome\cdot\Ome=0$. 
Thus $\b\cdot\ome$ is a section of $T^{1,0}(-D)\otimes\w^{0,1}$ which gives the class $[\b\cdot\ome]\in H^1(S, \Theta(-D))$. 
 
Then applying 
unobstructed deformations in the theorem \ref{th: unobstructedness  theorem}, we obtain the following, 
\bgn{theorem}\label{th: bihermitian on compact Kahler surfaces}
Let $S$ be a compact K\"ahler surface with complex structure $J$ and 
K\"ahler form $\ome$.
We denote by $g$ the K\"ahler metric on the K\"ahler surface $(S, J,\ome)$.
If there is a non-zero holomorphic Poisson structure $\b$ on $S$, 
then the surface $S$ admits deformations of bihermitian structures 
$(J, J^-_t, h_t)$ which satisfies 
$J^-_0=J$, $h_0=g$ and 
\bgn{equation}\label{eq: bihermitian on compact Kahler surfaces}
\frac{d}{dt}J^-_t|_{t=0} =-2(\b\cdot\ome+\ol\b\cdot\ome),
\end{equation}
where $\b\cdot \ome$ is the $\ol\pa$ closed $T^{1,0}$-valued forms of type $(0,1)$ which gives the Kodaira-Spencer class 
$-2[\b\cdot \ome]\in H^1(S, \Theta)$ of the deformations $\{J^-_t\}$.
In particular, if the class $[\b\cdot \ome]\in H^1(S, \Theta)$ does not vanish, then $(J, J^-_t, h_t)$ is a distinct bihermitian structure for small $t \neq 0$.
\end{theorem}
\bgn{proof}
 For a family of sections $a(t)$ of $\CL^2(-D)$ and real sections 
 $b(t)$ of $(\ol L^+_\J\cdot L^-_\J\oplus L^+_\J\cdot\ol L^-_\J)^\R$,
 we define a family of section $Z(t)$ of $\CL^2$ by $e^{Z(t)}= e^{a(t)}\, e^{b(t)}$, where 
 we denote by $(\ol L_\J^+\cdot L_\J^-\oplus L_\J^+\cdot\ol L_\J^-)^\R$ the real subbundle of the bundle $
 (\ol L_\J^+\cdot L_\J^-\oplus L_\J^+\cdot\ol L_\J^-)$.
 Since $b(t)\cdot \Ome =0$, 
 we have 
 $$
 de^{Z(t)}\cdot \Ome =d e^{a(t)}\, e^{b(t)}\cdot\Ome=de^{a(t)}\cdot \Ome
 $$
  Then since $a(t)$ is a section of $\CL^2(-D)$, it turns out that 
  $de^{Z(t)}\cdot \Ome$ is a smooth differential form on $S$.
  The action of $Z(t)\in \CL^2$ gives rise to almost bihermitian structures 
  $(J^+_t, J^-_t, h_t)$ with $J^\pm_0=J$ and $h_0=g$.
  
We shall construct $a(t)$ and $b(t)$ which satisfy the following three equations,
 \bgn{align}
 &d e^{Z(t)}\cdot\Ome =0\\
 &de^{Z(t)}\cdot\psi=0 \\
 &J^+_t=J,
 \end{align}
As in section 5, the structure $J^\pm_t$ is described by the adjoint action of a section $\Gam^\pm_t=\Gam^\pm(a(t), b(t))\in $GL$(T)$.
Then the equation $J^+=J$ is equivalent to $[\Gam^+_t, J]=[\Gam(a(t), b(t)), J]=0$.
We denote by $\(d e^{Z(t)}\)_{[i]}$ the $i$-th term of $\(d e^{Z(t)}\)$ in $t$ and 
$(\Gam^+_t)_{[i]}$ is also the $i$-th term of $(\Gam^+_t)$ in $t$.
Thus the three equations are reduced to the following equations for all integer $i\geq 0$: 
\bgn{align}
&\(d e^{Z(t)}\)_{[i]}\cdot\Ome =0\\
 &\(de^{Z(t)}\)_{[i]}\cdot\psi=0 \\
&[(\Gam^+_t)_{[i]}, J]=0.
\end{align}
We shall construct our solutions by the induction on $t$.

At first, we set $\h a_1:=\b+\ol\b$. 
Then the proposition \ref{prop: GK deformations} yields a real section $\h b_1\in \(\ol L_\J^+\cdot L_\J^-\oplus L_\J^+\cdot\ol L_\J^-\)^\R$ such that 
\bgn{align}
&d(\h a_1+\h b_1)\cdot\Ome =d\b \cdot \Ome =0\\
&d(\h a_1+\h b_1)\cdot\psi=0,
\end{align}
where we set $\b\cdot\Ome=1.$
We denote by $\Gam^+(\h a_1, \h b_1)$ the first term $(\Gam^+_t)_{[1]}$ in $t$ for $\h a_1,\h b_1$.
Then from the lemma \ref{lem: Gam-k}, we have  
\bgn{align}
\Gam^+(\h a_1, \h b_1) =\pi\circ \((\ad_{\h a_1}+\ad_{\h b_1})\circ 
\h J^+\circ \Ad_{e^\ome} \)
\end{align}
As in the lemma \ref{lem: ad gam k}, 
we define $\gam_1\in T\cdot T^*$ and $\ad_{\zeta_1}$ by 
\bgn{align}
\ad_{\zeta_1}:&=\[\, \Gam^+(\h a_1, \h b_1), J\, \]\\
\ad_{\gam_1}:&=\frac{-1}{(2\sqrt{-1})^2}\,\[\ad_{\zeta_1}, J\, \]
\end{align} 
Then it follows from the lemma \ref{lem: ad gam k} that 
$\gam_1$ is a section of $\CL^2(-D)$ and we have 
\bgn{align}\label{ali: ad gam1}
[\ad_{\gam_1}, J]+ \[\, \Gam^+(\h a_1, \h b_1), J\, \]=&
\frac{-1}{(2\sqrt{-1})^2}\[\[\ad_{\zeta_1}, J\, \],\, J\]+\ad_{\zeta_1}=0
\end{align}
From the lemma \ref{lem: rho k}, we have 
$d\gam_1\cdot\Ome =-d\rho_1\Ome$, where 
$\rho_1$ is a function with $\rho_1\Ome $ is a smooth form, that is, 
$\rho_1$ is a section of $\w^0\ol\L(-D)$.
Then we define $a_1$ by 
\bgn{align}
a_1=\h a_1+\gam_1+\rho_1
\end{align}
Then we have $da_1\Ome=0$.
Then applying the proposition \ref{prop: GK deformations} again, we have a section $b_1$ of 
$(\ol L^+\cdot L^-\oplus L^+\cdot\ol L^-)^\R$ such that $d(a_1+b_1)\cdot\psi=0$. 
From the lemma \ref{lem: Gam-k} and the lemma \ref{lem: gam}, we have 
\bgn{align}
\Gam^+( a_1, b_1)=&\pi\circ(\ad_{\h a_1}+\ad_{b_1}+\ad_{\gam_1})\circ\h J^+\circ\Ad_{e^\ome}\\
=&\pi\circ(\ad_{\h a_1}+\ad_{ b_1})\circ\h J^+\circ\Ad_{e^\ome}
+\ad_{\gam_1}
\end{align}
From the lemma \ref{lem: key b }, we have 
$$
[\pi\circ\ad_{b_1}\circ\h J^+\circ\Ad_{e^\ome}, J]=
[\pi\circ\ad_{\h b_1}\circ\h J^+\circ\Ad_{e^\ome}, J]=0
$$
Then it follows from (\ref{ali: ad gam1}) that
\bgn{align}
\[\, \Gam^+( a_1, b_1), J\, \]=&  
\[\, \Gam^+(\h a_1, \h b_1), J\, \]+
[\ad_{\gam_1}, J] =0,
\end{align}
Thus we obtain 
\bgn{align}
&d(a_1+b_1)\cdot\Ome=0\\
&d(a_1+b_1)\cdot\psi=0\\
&\[ \Gam^+(a_1, b_1), J]=0
\end{align}

Next we assume that there is a set of sections $a_1,\cdots, a_{k-1}$ of $\CL^2(-D)$ and sections $b_1,\cdots, b_{k-1}$ of $
(\ol L_\J^+\cdot L_\J^-\oplus L_\J^+\cdot\ol L_\J^-)^\R$ such that 
\bgn{align}
 &\(d e^{Z(t)}\)_{[i]}\cdot\Ome =0\\
 &\(de^{Z(t)}\)_{[i]}\cdot\psi=0 \\
 &\[\, \Gam^+(a(t), b(t))_{[i]}, J\]=0,
 \end{align}
 for all $0\leq i <k$, where $\(d e^{Z(t)}\)_{[i]}$ denotes the $i$-th term of 
 $\(d e^{Z(t)}\)$ in $t$ and $\Gam^+(a(t), b(t))_{[i]}$ is the $i$-th term of 
 $\Gam^+(a(t), b(t))$ for ${\displaystyle a(t) =\sum_{j=1}^{k-1} \frac{t^j}{j!}a_j, \,\, b(t)=\sum_{j=1}^{k-1} \frac{t^j}{j!}b_j}$. 
Then the proposition \ref{prop: GK deformations} yields a section $\h a_k$ of $\CL^2(-D)$ and a section $\h b_k$ of 
$(\ol L_\J^+\cdot L_\J^-\oplus L_\J^+\cdot\ol L_\J^-)^\R$ such that 
\bgn{align}
&\( de^{\h Z(t)}\cdot \Ome \)_{[k]} =\( de^{\h a(t)}\cdot\Ome \)_{[k]}=0,\\
&\(de^{\h Z(t)}\cdot\psi\)_{[k]} =0
\end{align}
where $\h Z(t)$ is a section of $\CL^2$ given by 
$e^{\h Z(t)}=e^{\h a(t)}\, e^{\h b(t)}$ and 
$$\h a(t) =\sum_{j=1}^{k-1} \frac{t^j}{j!} a_j + \frac{t^k}{k!} \h a_k
,\qquad \h b(t) =\sum_{j=1}^{k-1} \frac{t^j}{j!} b_j + \frac{t^k}{k!} \h b_k$$
Then for the section $\Gam^+(\h a(t),\h b(t))$ of GL$(T)$, 
as in the lemma \ref{lem: ad gam k} we define $\gam_k\in T\cdot T^*$ and $\ad_{\zeta_k}$ by 
\bgn{align}
&\ad_{\zeta_k}:=\[ \Gam^+(\h a(t),\h b(t))_{[k]}, \, J\]\\
&\ad_{\gam_k}:= \frac{-k!}{(2\sqrt{-1})^2}\[\ad_{\zeta_k}, \,J\]
\end{align}
Then from the lemma \ref{lem: ad gam k}, we see that $\gam_k$ is a section of $\CL^2(-D)$ and 
we have 
\bgn{align}\label{ali: [ad gam k, J]+}
&\frac1{k!}\[\ad_{\gam_k}, J\] +\[ \Gam^+(\h a(t),\h b(t))_{[k]}, \, J\]=0
\end{align}
The lemma \ref{lem: rho k} shows that $d\gam_k\cdot \Ome =-d\rho_k\cdot\Ome$ for 
a global function $\rho_k$ of $\w^0\ol\L(-D)$.
We define $a_k \in \CL^2(-D)$ by 
\bgn{align}\label{ali: ak definition}
a_k :=\h a_k+\gam_k+\rho_k
\end{align}
Then we have 
$$\(de^{a(t)}\)_{[k]}\Ome =\( de^{\h a(t)}\)_{[k]}\cdot\Ome +d( \gam_k+\rho_k)\cdot\Ome=0. 
$$ 
Applying the proposition \ref{prop: GK deformations} again,  we have a section 
$b_k$ of $\(\ol L_\J^+\cdot L_\J^-\oplus L_\J^+\cdot\ol L_\J^-\)^\R$ with 
$\( de^{Z(t)}\cdot \psi\)_{[k]}=0$, where $Z(t)=\log{\( e^{a(t)}\, e^{b(t)}\)}$. 
As in lemma \ref{lem: Gam-k}, $(\Gam^+_t)_{[k]}=\Gam^+(a(t), b(t))_{[k]}$ satisfies the following, 
\bgn{align}\label{ali: [Gam+(a(t), b(t))[k],J]}
[\Gam^+(a(t),b(t))_{[k]}, J] = &\frac 1{k!}[\pi\circ(\ad_{a_k}+\ad_{b_k})\circ\h J^\pm\circ\Ad_{e^\ome}, \, J]+[\wtil{\Gam^\pm_k}(a_{<k}, b_{<k}), \, J]
\end{align}
Substituting (\ref{ali: ak definition}) into (\ref{ali: [Gam+(a(t), b(t))[k],J]}) and using lemma \ref{lem: key b } and lemma \ref{lem: gam}, we have
\bgn{align*}
[\Gam^+(a(t),b(t))_{[k]}, J]=&
\frac 1{k!}[\pi\circ(\ad_{\h a_k}+\ad_{\gam_k})\circ\h J^\pm\circ\Ad_{e^\ome}, \, J]+[\wtil{\Gam^\pm_k}(a_{<k}, b_{<k}), \, J]\\
=&\frac 1{k!}[\pi\circ\ad_{\h a_k}\circ\h J^\pm\circ\Ad_{e^\ome}, \, J]+[\wtil{\Gam^\pm_k}(a_{<k}, b_{<k}), \, J] +\frac1{k!}[\ad_{\gam_k}, J]
\end{align*}
From lemma \ref{lem: Gam-k} and lemma \ref{lem: key b } we also have 
\bgn{align}
[\Gam^+(\h a(t),\h b(t))_{[k]}, J]=&\frac 1{k!}[\pi\circ\ad_{\h a_k}\circ\h J^\pm\circ\Ad_{e^\ome}, \, J]+[\wtil{\Gam^\pm_k}(a_{<k}, b_{<k}), \, J] 
\end{align}
Thus from(\ref{ali: [ad gam k, J]+}) we obtain
\bgn{align}
[\Gam^+(a(t),b(t))_{[k]}, J] = \frac1{k!}\[\ad_{\gam_k}, J\] +\[ \Gam^+(\h a(t),\h b(t))_{[k]}, \, J\]=0,
\end{align}
 where $a(t) =\sum_{j=1}^k \frac{t^j}{j!} a_j$ and 
 $b(t)=\sum_{j=1}^k\frac{t^j}{j!}b_j$. 
 Thus $Z(t)$ satisfies the equations, 
 \bgn{align}
&\(d e^{Z(t)}\)_{[k]}\cdot\Ome =0\\
&\(de^{Z(t)}\)_{[k]}\cdot\psi=0 \\
&\[\, \Gam(a(t), b(t))_{[k]}, J\]=0,
\end{align}
In section 6, we shall show that the formal power series $Z(t)$ is a convergent series which is smooth. 
Then the sections $a(t)$ and $b(t)$ give 
deformations of bihermitian structures 
$(J^+_t, J^-_t, h_t)$.
Finally we shall show that the family of deformations satisfies the equation 
(\ref{eq: bihermitian on compact Kahler surfaces}) in 
the theorem \ref {th: bihermitian on compact Kahler surfaces}. We already have $[\Gam^+_t, J]=0$ which implies that $J^+_t=J$. 
From the lemma \ref {lem: Gam-k} and the lemma \ref{lem: key b }, the 1st term of $J^-_t$ is given by 
\bgn{align*}
[(\Gam^-_t)_{[1]},J]=&[\pi\circ (\ad_{\h a_1} +\ad_{\gam_1} +\ad_{\h b_1})\circ \h J^-\circ\Ad_{e^\ome}), \, J]\\
=&[(\ad_{\gam_1}+\pi\circ\ad_{\h a_1}\circ\h J^-\circ\Ad_{e^\ome}),\, J]
\end{align*}
Since $\h a_1=\b+\ol\b$, we have $\pi \circ\ad_{\h a_1} |_T=0$. 
We also have  
$$[\ad_{\gam_1}, J] 
= [(\pi\circ \ad_{\gam_1}\circ J^*\circ \ad_\ome), J]=-[\Gam^+(\h a_1, \h b_1), J].$$
Thus we obtain 
$$
[(\Gam^-_t)_{[1]}, J]=2[(\pi\circ \ad_{\h a_1} \circ J^*\circ\ad_{\ome}),\, J]
$$
Then we have for a vector $v$,
\bgn{align*}
2(\pi\circ\ad_{\h a_1} \circ J^*\circ \ad_{\ome})v 
=&-2\[\b+\ol\b , \,[\ome, Jv]\]\\
=&-2\[[\b+\ol\b,\,\ome], \,Jv\]=-2(\b\cdot\ome+\ol\b\cdot\ome)Jv.
\end{align*}
Thus it follows that $\frac d{dt} J^-_t|_{t=0} =[(\Gam^-_t)_{[1]}, J]=
-2(\b\cdot\ome+\ol\b\cdot\ome)$ and the Kodaira-Spencer class of deformations $\{J^-_t\}$ is given by the class 
$-2[\b\cdot \ome]\in H^1(M.\Theta)$. 
If the class $[\b\cdot \ome]\in H^1(M.\Theta)$ does not vanish, then the deformations $\{J^-_t\}$ is not trivial. 
Thus $(X, J^-_t)$ is not biholomorphic to $(X, J)$ for small $t\neq 0$.

Hence we have the result.
\end{proof}
 \bgn{theorem} \label{th: bihermitian and Poisson}
A compact K\"ahler surface admits non-trivial bihermitian structure 
with the torsion condition and the same orientation if and only if $S$ has nonzero holomorphic Poisson structure.
 \end{theorem}
 \bgn{proof}
 It is already shown in \cite{A.G.G}, \cite{Hi2} that 
 a non-trivial bihermitian structure with the torsion condition and the same orientation carries a non-zero holomorphic Poisson structure. 
 It follows from the theorem \ref{th: bihermitian on compact Kahler surfaces} that if a compact K\"ahler surface $S$ has a non-zero holomorphic Poisson structure, then $S$ admits a nontrivial bihermitian structure. Thus the result follows.
 \end{proof}
 
\section{The convergence} 
In order to show the convergence of the power series in section 6, 
we apply the similar method in \cite{Go1}, \cite{Go3}.
We also use the same notation as in \cite{Ko}.
Let $P(t)=\sum_k P_k t^k$ be a power series in $t$ whose coefficients are
sections of a vector bundle on a Riemannian manifold.
We denote by $\|P_k\|_s$ the Sobolev norm of the section $P_k$ 
which is given by the sum of the $L^2$-norms of  $i$ th derivative of $P_k$ for all $i\leq s$, where $s$ is a positive integer with $s>2n+1$.
We put $\|P(t)\|_s=\sum_k \|P_k\|_st^k$.
Given two power series $P(t), Q(t)$, if 
$\|P_k\|\leq \|Q_k\|$ for all $k$, then we denote it by 
$$
P(t)\<Q(t).
$$
For a positive integer $k$, if $\|P_i\|\leq \|Q_i\|$ for all $i\leq k$,
we write it by 
$$
P(t)\underset k{\<} Q(t).
$$
We also use the following notation.
If $P_i=Q_i$ for all $i\leq k$, we write it by 
\bgn{equation}\label{eq: notation 1}
P(t)\underset k{\equiv} Q(t).
\end{equation}
Let $M(t)$ be a convergent power series defined by 
\bgn{equation}\label{eq: M(t)}
M(t) =\sum_{\nu=1}^\infty \frac1{16c}\frac{(ct)^\nu}{\nu^2}=\sum_{\nu=1}^\infty M_\nu t^\nu,
\end{equation}
for a positive constant $c$, which is determined later suitably.
The key point is the following inequality,
\bgn{equation}\label{eq: M(t)}
M(t)^2 \< \frac 1c M(t)
\end{equation}
We put $\lam =c^{-1}$. 
Then we also have 
\bgn{equation}\label{e M(t)}
e^{M(t)}\<\frac 1\lam e^\lam M(t).
\end{equation}
We will take $\lam$ sufficiently small which will be determined later.
(Note that $\lam$ gives a change of parameter $t$ by constant multiplication.)

As in the proposition \ref {prop: GK deformations}, 
if there is a set of sections $a_1,\cdots a_{k-1}$ of CL$^2$ which satisfies 
$$
\pi_{U^{-n+3}}\(e^{-a(t)}\, de^{a(t)}\)_{[i]}=0, \qquad \text{\rm for all }
i<k,
$$
and $\|a(t)\|_s\<_{k-1} K_1M(t)$,
then there is a set of real sections $b_1,\cdots, 
b_k\in (L^-_{\J}\cdot\ol L^+_{\J}\oplus \ol L^-_{\J}\cdot L^+_{\J})$    which satisfy 
the following equations:
\bgn{align}
&\pi_{U^{-n+3}}\(e^{-Z(t)}\, de^{Z(t)}\)_{[k]}=0\\
&\( d e^{Z(t)}\cdot\psi_0\)_{[i]}=0,\qquad \text{\rm for all } i\leq k\\
&\frac1{k!}\|\h a_k\|_s < K_1\lam M_k\\
&\frac1{k!}\|\h b_k\|_s <K_2M_k
\end{align}
where  $\h a_k$ is the section in the proposition \ref{prop: k th-G.complex} and
 $M(t)$ is the convergent series in (\ref{eq: M(t)}) with a constant $\lam$. 
 Note that 
$K_1$ is a positive constant and a positive constant $K_2$ is determined by $\lam, K_1$.
We also have an estimate of $e^{Z(t)}=e^{a(t)}e^{b(t)}$ in \cite{Go1}, 
$$
\|Z(t)\|\< _k M(t).
$$

Then $\gam_k$ in the lemma \ref{lem: ad gam k} satisfies 
\bgn{align}
\|\gam_k\|_s <& \| \Gam^+_k (a_{<k}, \h a_k, b_{<k}, \h b_k) \|_s\\
<&2 \|\h a_k\|_s +2\|\h b_k\|_s+\| \wtil{\Gam^+_k }(a_{<k}, b_{<k})\|_s
\end{align}
Recall that $\Gam_t^+ =\pi \( \Ad_{e^{Z(t)}}\circ \h J^+\circ\Ad_{e^\ome}\)$.
Then we have an estimate of the non-linear term $\| \wtil{\Gam^+_k }(a_{<k}, b_{<k})\|_s$ 
$$\| \wtil{\Gam^+_k }(a_{<k}, b_{<k})\|_s < C\|( e^{Z(t)}-Z(t) -1)_{[k]}\|_s,$$
where $C$ denotes a constant. 
It follows from (\ref{e M(t)}) that $\|( e^{Z(t)}-Z(t) -1)_{[k]}\|_s< C(\lam)M_k$, where $C(\lam)$ satisfies 
$\lim_{\lam\to0}C(\lam)=0$. 
Thus we have 
$$
\frac1{k!}\|\gam_k\|_s< \frac{2}{k!}\(\|\h a_k\|_s+\|\h b_k\|_s\)+ C(\lam)M_k < 2\lam K_1 M_k +
2K_2M_k+C(\lam)M_k.
$$
By using the Hodge decomposition and the Green operator
We also have a unique global function $\rho_k$ of $\w^0\ol\L(-D)$ which satisfies the followings, 
\bgn{align}
&d\rho_k\Ome =-d\gam_k\cdot\Ome\\
&\|\rho_k\|_s\leq C_1\|\gam_k\|_s,
\end{align}
where $C_1$ is a constant.
Then we obtain 
\bgn{align*}
\frac1{k!}\|a_k\|_s < &\frac1{k!}\|\h a_k \|_s+ \frac1{k!}\|\gam_k\|_s
+\frac1{k!}\|\rho_k\|_s\\
<&\frac1{k!}\|\h a_k \|_s+ \frac1{k!}(1+C_1)\|\gam_k\|_s\\
<&\lam K_1 M_k +2(1+C_1)(\lam K_1 M_k +K_2M_k+C(\lam)M_k)
\end{align*}
We take $\lam $ and $K_2$ sufficiently small such that $\lam K_1 M_k +2(1+C_1)(\lam K_1 M_k +K_2M_k+C(\lam)M_k)< K_1M_k$. 
Then we obtain
$$\frac1{k!}\|a_k\|_s < \frac1{k!}\|\h a_k \|_s+ \frac1{k!}\|\gam_k\|_s< K_1 M_k.$$
Thus our solution $a(t)$ satisfies that 
$\|a(t)\|_s \<_k K_1M(t)$ for all $k$ by the induction.
It implies that $a(t)$ is a convergent series. 
Applying the proposition \ref{prop: GK deformations} again, we have 
$\|b(t)\|_s<_k K_2 M(t)$. 
Hence $b(t)$ is also a convergent series.
Thus it follows that $Z(t)$ is a convergent series.

\section{Applications}
\subsection{Bihermitian structures on del Pezzo surfaces}
A del Pezzo surface is by definition a smooth algebraic surface with ample anti-canonical line bundle.
A classification of del Pezzo surfaces are well known, they are 
$\C P^1\times \C P^1$ or $\C P^2$ or a surface $S_n$ which is the blow-up of $\C P^2$ at $n$ points $P_1,\cdots, P_n$, $(0< n \leq 8)$. 
The set of the points $\Sig:=\{P_1,\cdots, P_n\}$ must be in
{\it general position} to yield a del Pezzo surface.
The following theorem is due to Demazure, {\rm \cite{De} (see page 27)}, which shows the meaning of {\it general position},
\bgn{theorem}
The following conditions are equivalent:\\
(1) The anti-canonical line bundle of $S_n$ is ample \\
(2) No three of $\Sig$ lie on a line, no six of  $\Sig$ lie on a conic and 
no eight of $\Sig$ lie on a cubic with a double point $P_i\in \Sig$\\
(3) There is no curve $C$ on $S_n$ with $-K_{S_n}\cdot C \leq 0$.\\
(4) There is no curve $C$ with $C\cdot C=-2$ and $K_{S_n}\cdot C=0$.
\end{theorem}
\bgn{remark}
If three points lie on a line $l$, then the strict transform $\h l$ of $l$ in $S_3$ is a $(-2)$-curve with $K_{S_3}\cdot \h l=0$. 
If six points belong to a conic curve $C$, then 
the strict transform form $\h C$ of $C$ is again a $(-2)$-curve with $K_{S_6}\cdot \ol C=0$.
If eight points $P_1\cdots, P_8$ lie on a cubic curve with a double point $P_1$, then the strict transform $\h C$ of $C$ satisfies
$\h C \thicksim\pi^{-1}C -2E_1-E_2-\cdots -E_8$, where 
$E_i$ is the exceptional curve $\pi^{-1}(P_i)$. 
Then we also have $\h C^2=-2$ and $K_{S_8}\cdot \h C=0$.
\end{remark}

Let $D$ be a smooth anti-canonical divisor of ${S_n}$ which
 is given by the zero locus of a section $\b\in H^0({S_n}, K_{S_n}^{-1})$. 
Since the anti-canonical bundle $K_{S_n}^{-1}$ is regarded as the bundle of $2$-vectors $\w^2\Theta$ 
and $[\b, \b]_{S}=0\in \w^3\Theta$ on $S_n$,  every section $\b$ is  a holomorphic Poisson structure. 
On $S_n$, we have the followings,
$$
\dim H^1(S_n, \Theta)=
\bgn{cases}
&2n-8\quad ( n=5, 6,7,8) \\
&0\quad \,\,\,\qquad (n<5)
\end{cases}
$$
$$
\dim H^0(S_n, K^{-1}) =10-n
$$
and 
$$
H^{1,1}(S_n ) =1+n.
$$

Further we have 
$H^2(S_n ,\Theta)=\{0\}, \quad H^1(S_n , \w^2\Theta)\cong 
H^1(S_n , -K_{S_n})=\{0\}$.
Hence the obstruction vanishes and we have deformations of \complex structures parametrized by $H^0(S_n ,K_{S_n}^{-1})\oplus H^1(S_n , \Theta)$.

In particular,  if $n\geq 5$, we have deformations of ordinary complex structures on $S_n$
\bgn{proposition}\label{prop: del Pezzo 1}
Let $D$ be a smooth anti-canonical divisor given by the zero locus of $\b$ as above. 
Then there is a K\"aher form $\ome$ with the class 
$[\b\cdot \ome ]\neq 0\in H^1(S_n, \Theta)$.
\end{proposition}

We also have $H^2(\C P^1\times \C P^1,\Theta)=0$ and 
$H^1(\C P^1\times \C P^1, -K)=0$.

Thus we can apply our construction to every del Pezzo surface.
From the main theorem together with the proposition
\ref{prop: del Pezzo 1}, we have 
\bgn{proposition}
Every del Pezzo surface admits deformations of bihermitian structures 
$(J, J^-_t, h_t)$ with $J^-_0=J$ which satisfies 
\bgn{equation}
\frac d{dt} J^-_t|_{t=0} =-2(\b\cdot\ome+\ol\b\cdot\ome),
\end{equation}
for every K\"ahler form $\ome$ and every holomorphic Poisson structure $\b$.
Further, a del Pezzo surface $S_n$ $(n\geq 5)$ admits 
distinct bihermitian structures $(J, J^-_t, h_t)$, that is, 
 the complex manifold $(X, J^-_t)$ is not biholomorphic to $(X, J)$ for small $t\neq 0 $.
\end{proposition}
Note that for small $t\neq 0$, $J^-_t\neq \pm J$.
We will give a proof of the proposition \ref{prop: del Pezzo 1} in the rest of this subsection. 

Let $N_D$ is the normal bundle to $D$ in ${S_n}$ and $i^*T_{S_n}$ the pull back of the tangent bundle $T_{S_n}$ of ${S_n}$ by the inclusion $i : D\to {S_n}$. 
Then we have the short exact sequence, 
$$
0\to T_D \to i^*T_{S_n} \to N_D\to 0
$$
and we have the long exact sequence
\bgn{align*}
0\to &H^0(D,T_D)\to 
H^0(D, i^*T_{S_n})\to 
H^0(D, N_D)\overset{\delta}\to H^1(D, T_D)\to \cdots \\
\end{align*}
Since the line bundle $N_D$ is positive, $H^1(D, N_D)=\{0\}$ and 
$\dim H^0(D, N_D)$ is equal to the intersection number $D\cdot D=9-n$ by the Riemann-Roch theorem.
Since $D$ is an elliptic curve, 
$\dim H^1(D, T_D) =\dim H^0(D, T_D) =1$. 
Hence if follows that 
\bgn{equation}\label{eq: key inequ in del Pezzo}
9-n\leq \dim H^0(D, i^*T_{S_n})\leq 10-n.
\end{equation}
Let $\I_D$ be the ideal sheaf of $D$ and $\O_D$ the structure sheaf of $D$. 
Then we have the short exact sequence
$$
0\to \I_D\to \O_{S_n}\to i_*\O_D\to 0
$$
By the tensor product, we also have 
\bgn{equation}\label{eq: exact sequence 2 in del Pezzo}
0\to \I_D\otimes T_{S_n}\to T_{S_n}\to i_*\O_D\otimes T_{S_n}\to 0
\end{equation}
Then from the projection formula we have 
$$
H^p({S_n}, i_*\O_D\otimes T_{S_n})\cong H^p({S_n}, i_*(\O_D\otimes i^*T_{S_n}))
\cong H^p(D, i^*T_{S_n}),
$$
for $p=0,1,2$.
From (\ref{eq: exact sequence 2 in del Pezzo}), we have the long exact sequence,
\bgn{align}
H^0({S_n}, T_{S_n})\to H^0(D, i^*T_{S_n})\to H^1({S_n}, \I_D\otimes T_{S_n})
\overset{j}\to H^1({S_n}, T_{S_n})\to\cdots
\end{align}
Hence we obtain
\bgn{lemma}\label{lem: key lemma in del Pezzo}
 The map $j : H^1({S_n}, \I_D\otimes T_{S_n})\to H^1({S_n}, T_{S_n})$
is not the zero map .
\end{lemma}
\bgn{proof}
We have the exact sequence,
\bgn{align}
\cdots\to H^0(D, i^*T_{S_n})\to H^1({S_n}, \I_D\otimes T_{S_n})\overset{j}\to H^1({S_n}, T_{S_n})
\end{align}
From the Serre duality with $\I_D=K_{S_n}$, we have 
$H^0(S_n, \I_D\otimes T_{S_n})\cong H^2(S_n, \Ome^1_{S_n})=\{0\}$ and 
$H^2(S_n, \I_D\otimes T_{S_n})=H^0(S_n, \Ome^1)=0$. 
From the Riemann-Roch theorem, 
$\dim H^1(S_n, \I_D\otimes T_{S_n})=n+1$.
Then it follows from (\ref{eq: key inequ in del Pezzo}) that 
$$
\dim H^0(D , i^*T_{S_n})<\dim H^1(S_n, \I_D\otimes T_{S_n})
$$
Note $10-n<n+1$ for all $n\geq 5$.
Hence the map $j$ is non zero.
\end{proof}
\bgn{remark}
Since $n\geq 5$, we have 
$H^0({S_n}, T_{S_n})=\{0\}$.
Applying the Serre duality with $K_{S_n}=\I_D$, we have $H^2(S_n, T_{S_n})\cong H^0(S_n, \I_D\otimes\Ome^1)=0$.
From the Riemann-Roch, we obtain $\dim H^1({S_n}, T_{S_n})=2n-8$.
\end{remark}
Let $\b$ be a non-zero holomorphic Poisson structure $S_n$ with the smooth divisor $D$ as the zero locus. 
Then $\b$ is regarded as a section of $\I_D\otimes \w^2\Theta$. 
Thus the section $\b\in H^0({S_n}, \I_D\otimes\w^2\Theta)$ gives an identification, 
$$\Ome^1\cong \I_D\otimes T_{S_n}.$$
Then the identification induces the isomorphism 
$$\h \b : H^1({S_n}, \Ome^1) \cong H^1({S_n}, \I_D\otimes T_{S_n}).$$
Let $j$ be the map in the lemma \ref{lem: key lemma in del Pezzo}.
Then we have the composite map $j\circ \h\b :H^1(S_n ,\Ome^1)\to H^1(S_n, \Theta)$ 
which is given by 
the class $[\b\cdot \ome] \in H^1({S_n}, T_{S_n})$ for 
$[\ome] \in H^1({S_n}, \Ome^1)$.
\bgn{proposition}
The composite map $j\circ \h \b : H^1({S_n}, \Ome^1)\to H^1({S_n}, T_{S_n})$ is not the zero map. 
\end{proposition}
\bgn{proof}
Since the map $\h\b$ is an isomorphism, $\h\b(\ome)$ is not zero.
It follows from lemma \ref{lem: key lemma in del Pezzo} that the map $j$ is non-zero.
Hence the composite map $j\circ \h \b$ is non-zero also.
\end{proof}
\bgn{proof} of lemma \ref{prop: del Pezzo 1}
The set of K\"ahler class is an open cone in 
$H^{1,1}({S_n}, \R)\cong H^2(S_n \R)$. We have the non-zero map 
$j\circ\h\b: H^{2}({S_n},\C) \cong H^1({S_n},\Ome^1) \to H^1({S_n},\Theta)$
for each $\b\in H^0({S_n}, K^{-1})$ with $\{\b=0\} =D$.
It follows that the kernel $j\circ\h\b$ is a closed subspace and 
 the intersection $\ker(j\circ\h\b)\cap H^{2}({S_n},\R)$ is closed in 
 $H^{2}({S_n}, \R)$ whose dimension is strictly less than $\dim H^{2}({S_n}, \R)$.
 Thus the complement in the K\"ahler cone 
 $$\big\{ [\ome] :\text{\rm K\"ahler class} \, |\, j\circ\h\b([\ome])\neq0\big\}$$
 is not empty. 
 Thus there is a K\"ahler form $\ome$ such that the class $[\b\cdot\ome]\in H^1(S_n ,\Theta)$ does not vanish for $n\geq 5$.
\end{proof}
We also remark that our proof of the lemma \ref {prop: del Pezzo 1} still works for degenerate del Pezzo surfaces. 
  \subsection{Vanishing theorems on surfaces}
 Let $M$ be a compact complex surface with canonical line bundle $K_M$.
  We shall give some vanishing theorems of the cohomology groups
  $H^1(M, -K_M)$ and $H^2(M, \Theta)$ on a compact smooth complex surface $M$, which are the obstruction spaces to deformations of \complex structures starting from the ordinary one $(X, \J_J)$.
The following is  practical to show the vanishing of $H^1(M, -K_M)$. 
 \bgn{proposition}\label{vanishing H1(-K)}
  Let $M$ be a compact complex surface with $H^1(M, {\cal O}_M)=0$.
  If $-K_M=m[D]$ for a irreducible, smooth curve $D$ with positive self-intersection number $D\cdot D>0$  and a positive integer $m$, then 
  $H^1(M, K^n_M)=0$ for all integer $n$.
  \end{proposition}
  The proposition is often used in the complex geometry. For completeness, we give a proof.
  \bgn{proof}
  Let $I_D$ be the ideal sheaf of the curve $D$. Then we have the short exact sequence, $0\to I_D\to {\cal O}_M\to j_*{\cal O}_D\to 0$, where $j: D\to X$.
 Then we have the exact sequence, 
 $$
 H^0(M, {\cal O}_M)\to H^0(M, j_*{\cal O}_D)\overset{\del}{\to}H^1(M, I_D)\to 
 H^1(M, {\cal O}_M) 
 $$
 It follows that the coboundary map $\del$ is a $0$-map. Thus  from $H^1(M, {\cal O}_M)=0$, we have $H^1(M,I_D)=H^1(M, -[D])=0$.
 We use the induction on $k$. We assume that $H^1(M, I_D^k)=H^1(M, -k[D])=0$ for a positive integer $k$.
 The short exact sequence $0\to I^{k+1}_D\to I^k_D\to j_*{\cal O}_D\otimes I^k_D\to 0$ induces the exact sequence, 
 $$
 H^0(M, j_*{\cal O}_D\otimes I_D^k)\to H^1(M, I_D^{k+1})\to H^1(M, I_D^k).
 $$
 By the projection formula, we have $H^0(M, j_*{\cal O}_D\otimes I_D^k)=H^0(D, -k[D]|_D)$. Since $D\cdot D>0$, it follows that the line bundle $ -k[D]|_D$ is negative and then
 $H^0(D,  -k[D]|_D)=H^0(M, I_D^k)=0$. It implies that $H^1(M, I_D^{k+1})=H^1(M, -(k+1)[D])=0$. 
 Thus by the induction, we have $H^1(M, -nD)=0$ for all positive integer $n$. 
 Applying the Serre duality, we have $H^1(M, -nD)\cong H^1(M, (n-m)D)=0$. 
 Thus $H^1(M, nD)=0$ for all integer $n$. 
 Then the result follows since $H^1(M, K^n)=H^1(M, -(nm)D)=0$.
  \end{proof}
 The author also refer to the standard vanishing theorem.
 If $D=\sum_i a_i D_i$ is a $\Q$-divisor on $M$, where $D_i$ is a prime divisor and $a_i \in \Q$. Let $\lceil a_i\rceil$ be the round-up of $a_i$ and $\lfloor a_i\rfloor$ the round-down of $a_i$. 
 Then the fractional part $\{a_i\}$ is $a_i -\lfloor a_i\rfloor$. 
 Then the round-up and the round-down of $D$ is defined by 
 $$
 \lceil D\rceil =\sum_i \lceil a_i\rceil D_i, \quad 
 \lfloor D\rfloor=\sum_i \lfloor a_i \rfloor D_i
 $$
and  $\{ D\}=\sum_i\{a_i\}D_i$ is the fractional part of $D$. 
A divisor $D$ is {\bf nef} if one has $D\cdot C\geq 0$ for any curve $C$. 
A divisor $D$ is {\bf nef} and {\bf big} if in addition, one has $D^2>0$.
We shall use the following vanishing theorem. 
The two dimensional case is due to Miyaoka and the higher dimensional cases
are due to Kawamata and Viehweg
\bgn{theorem} \label{th: vanishing K&V}
Let $M$ be a smooth projective surface and $D$ a $\Q$-divisor on $M$ such that \\
(1) supp$\{D\}$ is a divisor with normal crossings, \\
(2) $D$ is nef and big.\\
Then $H^i(M, K_M+\lceil D\rceil) =0$ for all $i>0$.
\end{theorem}
If $-K_M=mD$ is nef and big divisor where $D$ is smooth for $m>0$.
Then applying the theorem, we have 
$$
H^i(M, -K_M) \cong H^i(M, K_M-2K_M)=0,
$$
for all $i>0$.
\\
Next we consider the vanishing of the cohomology group $H^2(M, \Theta)$. 
Applying the Serre duality theorem, 
we have 
$$
H^2(M, \Theta)\cong H^0(M, \Ome ^1\otimes K_M)
$$
If $-K_M$ is an effective divisor $[D]$, then $K_M$ is given by the ideal sheaf 
$I_D$ of $D$. 
The short exact sequence: $0\to \Ome^1\otimes I_D \to \Ome^1\to 
\Ome^1\otimes{\cal O}_D\to 0$ gives us the injective map, 
$$
0\to H^0(M, \Ome^1\otimes K_M)\to H^0(M, \Ome^1).
$$ 
Hence we have 
\bgn{proposition}\label{vanishing H2Theta}
if $M$ is a smooth surface with effective anti-canonical divisor satisfying
$H^0(M, \Ome^1)=0$, then we have the vanishing $H^2(M, \Theta)=0$.
\end{proposition}

 \subsection{Non-vanishing theorem}
\bgn{proposition}\label{non-vanishing}
 Let $M$ be a K\"ahler surface with a K\"ahler form $\ome$ and 
 a non-zero Poisson structure $\b \in H^0(M, \w^2\Theta)$. 
 Let $D$ be the divisor defined by the section $\b$. 
 If there is a curve $C$ of $M$ with $C\,\cap$ supp $D=\emptyset$, 
 then the class $[\b\cdot\ome]\in H^1(M, \Theta)$ does not vanish.
 \end{proposition}
 \bgn{proof}
 Since $\b$ is not zero on the complement $M\bsh D$, there is a holomorphic symplectic form $\h \b$ on the complement.
 The symplectic form $\h \b$ gives the isomorphism 
$\Theta \cong \Ome^1$ on $M\bsh D$ which induces the isomorphism 
between cohomology groups
$H^1(M\bsh D, \Theta)\cong H^1(M\bsh D, \Ome^1)$.
Then the restricted class $[\b\cdot\ome] |_{M\bsh D}$ corresponds to the 
K\"ahler class $[\ome]|_{M\bsh D}\in H^1(M\bsh D, \Ome^1)
\cong H^{1,1}(M\bsh D)$ under the isomorphism. 
Since there is the curve $C$ on the complement $M\bsh D$ and $\ome$ is a K\"ahler form, the class $[\ome|_C]\in H^{1,1}(C)$ does not vanish.
Then it follows that the class $[\ome]|_{M\bsh D}\in H^1(M\bsh D, \Ome^1)$ does not vanish. 
It implies that $[\b\cdot\ome]|_{M\bsh D}$ does not vanish also.
Thus we have that the class $[\b\cdot \ome]\in H^1(M, \Theta)$ does not vanish.
 \end{proof}
\subsection{Deformations of bihermitian structures on the Hirzebruch surfaces $F_e$} 
Let $F_2$ be the projective space bundle of 
$T^*\C P^1\oplus {\cal O}_{\C P^1}$ over $\C P^1$ with fibre $\C P^1$,
$$
F_2=\Bbb P(T^*\C P^1\oplus {\cal O}_{\C P^1}).
$$

We denote by $E^+$ and $E^-$ the sections of $F_2$ with positive and negative self-intersection numbers respectively. 
An anti-canonical divisor of $F_2$ is given by $2E^+$, 
while the section $E^-$ with  
 $E^-\cdot E^-=-2$ is the curve which satisfies $E^+\cap E^-=\emptyset$.
 Thus we have the non-vanishing class $[\b\cdot \ome]\in H^1(F_2, \Theta)$, where $\b$ is a section of $-K$ with the divisor $2E^+$.
  \
(Note that the canonical holomorphic symplectic form $\h \b$ on the cotangent bundle $T^*\C P^1$ which induces the holomorphic Poisson structure $\b$. 
The structure $\b$ can be extended to $F_2$ which gives the anti-canonical divisor $2[E^+]$.)
\bgn{proposition}\label{prop: non-vanishing on F2}
The class $[\b\cdot\ome] \in H^1(F_2, \Theta)$ does not vanish for every K\"ahler form $\ome $ on $F_2$.
\end{proposition}
\bgn{proof}
The result follows from the proposition \ref{non-vanishing}.
\end{proof}

On the surface $F_2$, 
 the anti-canonical line bundle of $F_2$ is 
$2E^+$ and $H^1(F_2,{\cal O}_{F_2})=0$.
Hence from the proposition \ref{vanishing H1(-K)}, we  have the vanishing $H^i(F_2, -K_X)=\{0\}$ for all 
$i>0$. 
Since the surface $F_2$ is simply connected, it follows from the proposition 
\ref{vanishing H2Theta} that
 $H^2(F_2, \Theta)=0$. 
Hence the obstruction vanishes. 
It is known that every non-trivial small deformation of $F_2$ is $\C P^1\times \C P^1$. Thus we have
\bgn{proposition}
Let $(X, J)$ be the Hirzebruch surface $F_2$ as above.
Then there is a family of deformations of bihermitian structures $(J^+_t, J^-_t, h_t)$ with the torsion condition and the same orientation
with  $J^+_t=J^-_0=J$ such that 
$(X, J^-_t)$ is $\C P^1\times \C P^1$ for small $t\neq 0$.
\end{proposition}

Let $F_e$ be the projective space bundle  
 $\Bbb P({\cal O}\oplus {\cal O}(-e))$ over $\C P^1$
 with $e>0$. 
 There is a section $b$ with $b^2=-e$, which is unique if $e>0$. 
 Let $f$ be a fibre of $F_e$. 
 Then $-K$ is given by $2b+(e+2)f$, which is an effective divisor.
Thus from the proposition \ref{vanishing H2Theta}, 
we have $H^2(F_e, \Theta)=\{0\}$. 
$P_{-1}(F_e)=\dim H^0(F_e, K^{-1})$ is listed in the table 7.1.1 of \cite{Sa1}, 
$$
P^{-1}(F_e) =
\bgn{cases}
&9 \qquad\,\,\,\, e=0,1 \\
&9 \qquad\,\,\,\, e=2\\
&e+6\quad e\geq 3
\end{cases}
$$
Since $K$ is given by the ideal sheaf $I_D$ for the effective divisor $D=2b+(e+2)f$,
It follows from the Serre duality that $H^2(F_e, K^{-1}) =H^0(F_e, I_D^2)=\{0\}$.
Thus applying the Riemann-Roch theorem, we obtain
$$
\dim H^1(F_e, K^{-1}) =e -3,
$$
for $e\geq 3$. 
In the case $e=3$, note that $H^1(F_3, K^{-1})=H^2(F_3, \Theta)=\{0\}$.  From our main theorem, we have 
\bgn{proposition}
The Hirzebruch surface $F_e$ admits deformations of non-trivial bihermitian structures with the torsion condition and the same orientation
$(J, J^-_t, h_t)$ with $J^-_t\neq \pm J$ for small $t\neq 0$.
\end{proposition}
\subsection{Bihermitian structures on ruled surfaces 
$\Bbb P(T^*\Sig_g\oplus {\cal O}_{\Sig_g})$}

We can generalized our discussion of $F_2$ to
 the projective space bundle of $T^*M\oplus{\cal O}_M$
 over a compact K\"ahler manifold $M$. 
 Then we also have the Poisson structure $\b$  and as in the proposition \ref{non-vanishing}, it is shown that the class $[\b\cdot \ome]$
 does not vanish. 
 Thus we have the deformations of bihermitian structures from 
 the stability theorem \cite{Go2}.
  
 If $M$ is a Riemannian surface $\Sig_g$ of genus $g\geq 1$, then 
 the projective space bundle is called a ruled surface of degree $g$. 
 It is known that small deformations of any ruled surface of degree $g \geq 1$ 
 remain to be ruled surfaces of the same degree. 
We denote by $S=(X,J)$ the projective space bundle 
$\P(T^*\Sig\oplus {\cal O}_\Sig)$, where $X$ is the underlying differential manifold and $J$ is the complex structure.
 Applying our main theorem, we have
\bgn{theorem}\label{th: distinct bihermitian}
There is family of distinct bihermitian structures $(J, J^-_t, h_t)$ with the torsion condition and the same orientation on 
$S:=\P(T^*\Sig\oplus {\cal O}_\Sig)$, 
that is, the complex manifold $(X, J^-_t)$ is not biholomorphic to $S=(X, J)$ for small $t\neq 0$.
\end{theorem}

\subsection{Bihermitian structures on degenerate del Pezzo surfaces}
We shall consider the blow-up of $\C P^2$ at $r$ points which are not in 
general position. 
We follow the construction as in 
\cite {De}, (see page 36). 
 We have a finite set $\Sig=\{x_1,\cdots, x_r\}$ and $X(\Sig)$ obtained by successive blowing up at $\Sig$,
$$
X(\Sig)\to X(\Sig_{r-1})\to \cdots \to X(\Sig_1)\to \C P^2,
$$
At first $X(\Sig_1)$ is the blow-up of $\C P^2$ at a point $x_1\in \C P^2$
 and  we have $\Sig_i =\{x_1, \cdots, x_{i}\}$ and 
$X(\Sig_{i+1})$ is the blow-up of $X(\Sig_{i})$ at $x_{i+1}\in X(\Sig_i)$.
Let $E_i$ be the divisor given by  the inverse image of $x_i\in X(\Sig_{i-1})$.
If $\Gam$ is an effective divisor on $\C P^2$, one notes that mult$(x_i, \Gam)$ 
the multiplicity of $x_i$ on the proper transform of $\Gam$ in $X(\Sig_{i-1})$, 
and one says that $\Gam$ passes through $x_i$ if mult$(x_i, \Gam)>0$. 
 Define $\h E_1, \cdots, \h E_r$ by recurrence as follows, 
 On $X(\Sig_1)$, one put $\h E_1=E_1$ ; 
 on $X(\Sig_2)$, $\h E_1$ is a proper transform of the previous $E_1$ and 
 one also put $\h E_2=E_2$; on $X(\Sig_3)$, $\h E_1$ and $\h E_2$ are the proper transform of previous $\h E_1$ and $\h E_2$  respectively and $\h E_3=E_3$. 
 Then $\h E_1, \cdots, \h E_r$ are irreducible components of $E_1+\cdots+ E_r$. 
 
We assume that the following condition on $\Sig$, \medskip\\
(*)  For each $i=1,\cdots, r$, a point $x_i\in X(\Sig_{i-1})$ does not belong to a irreducible curve $\h E_j$ with self-intersection number $-2$ for $1\leq j\leq i-1$.\\ \\
If a point $x_i\in X(\Sig_{i-1})$ belongs to a irreducible curve 
$\h E_j$ with self-intersection number $-2$, then the proper
transform of $\h E_j$ becomes a curve with self-intersection number $-3$.
If there is a rational curve with self-intersection number $-3$ or less, the anti-canonical divisor of $X(\Sig)$ is not nef.
 
\bgn{definition} 
 A set of points $\Sig $ is in {\it almost general position} if $\Sig$ satisfies the following:\\
 (1) $\Sig$ satisfies the condition (*)   \\
 (2) No line passes through $4$ points of $\Sig$ \\
 (3) No conic passes through $7$ points of $\Sig$ \\
 \end{definition}
 We call $X(\Sig)$ {\it a degenerate del Pezzo surface} if $\Sig$ is  in almost general position.
Note that if $\Sig$ is in general position, $\Sig$ is in almost  general position. 
In \cite{De}, the following theorem was shown,\\
\bgn{theorem}{\rm \cite{De} }The following conditions are equivalent: \\
(1) $\Sig$ is in almost general position \\
(2) The anti-canonical class of $X(\Sig)$ contains a smooth and irreducible curve $D$. \\
(3) There is a smooth curve of $\C P^2$ passing all points of $\Sig$.\\
(4) $H^1(X(\Sig), K_{X(\Sig)}^n)=\{0\}$ for all integer $n$\\
(5) $-K_{X_{\Sig}}\cdot  C\geq 0$ for all effective curve $C$ on $X(\Sig)$ and in addition,
 if $-K_{X(\Sig)}\cdot C =0$, then $C\cdot C=-2$.
\end{theorem}
Then from (2) there is a smooth anti-canonical divisor on a degenerate del Pezzo surface and we have $H^1(X(\Sig),{\cal O}_X)=0$.
Hence from the proposition \ref{vanishing H1(-K)}, we have the vanishing $H^i(X(\Sig), -K_X)=0$, for all $i>0$. 
  A degenerate del Pezzo surface $X(\Sig)$ satisfies $H^0(X(\Sig),\Ome^1)=0$. 
 Then it follows from the proposition \ref{vanishing H2Theta} that $H^2(X(\Sig), \Theta)=0$.
 
Let $X(\Sig)$ be a degenerate del Pezzo surface which is not a del Pezzo surface, that is,  the anti-canonical class of $X(\Sig)$ is not ample. 
Then from (5), there is a $(-2)$-curve $C$ with $K_{X(\Sig)}\cdot C=0$. 
Then it follows that $C$ is a $\C P^1$.
Thus we contract $(-2)$-curves on a degenerate del Pezzo to obtain a complex surface with rational double points, which is called the Gorenstein log del Pezzo surface. 
Let $\b$ be a section of $-K_{X(\Sig)}$ with the smooth divisor $D$ as the zero set. 
We denote by $J$ the complex structure of the del Pezzo surface $X(\Sig)$.
From our main theorem, we have 
\bgn{theorem}
A degenerate del Pezzo surface admits deformations of distinct bihermitian structures $(J, J^-_t, h_t)$ with $J_0^-=J$ and $J^-_t\neq \pm J$ for 
small $t\neq 0$, that is, 
 $\frac d{dt} J^-_t|_{t=0}=
-2(\b\cdot \ome+\ol\b\cdot\ome)
$,
and 
the complex structure $J^-_t$ is not equivalent to 
$J$ of $X(\Sig)$ under diffeomorphisms for small $t\neq 0$, 
where $\ome$ is a K\"ahler form.
\end{theorem}
 \bgn{proof}
 If $X(\Sig)$ is a del Pezzo surface, we already have the result.
 If $X(\Sig)$ is not a del Pezzo but a degenerate del Pezzo, we still have $H^2(X(\Sig), \Theta)=
 H^1(X(\Sig), K^{-1})=\{0\}$. Thus we have deformations of bihermitian structures as in our main theorem. 
 It is sufficient to show that the class $[\b\cdot\ome]$ does not vanish.
 Since there is a smooth $(-2)$-curve $C$ with $K\cdot C=0$, the line bundle $K^{-1}|_C\to C\cong \C P^1$ is trivial. 
 If there is a point $P\in D\cap C$, then $\b(P)=0$ and it follows that $\b|_C\equiv 0$. Since the anti-canonical divisor $D$ is smooth, we have $D=C$. 
 However $D\cdot D=9-r\neq 0$ and $D\cdot C= -K\cdot C=0$. 
 It is a contradiction.
 Thus $D\cap C=\emptyset$. 
 Then applying the proposition \ref{non-vanishing}, 
 we obtain $[\b\cdot \ome]\neq 0\in H^1(X(\Sig), \Theta)$.
 \end{proof} 
\bgn{thebibliography}{99}
\bibitem{A.D}
V.~Apostolov an G.~Dloussky, 
{\it Bihermitian metrics on Hopf surfaces}, 
Math. Res Lett. 15 (2008) no. 5, 827-839
\bibitem{AN}
V.~Alexeev and V.~ Nikulin,
{\it Del Pezzo surfaces and K3 surfaces}, 
MSJ Memoirs, Mathematical Society of Japan, Vol. 15 (2006), MR 2227002,
Zbl 1097.1400, 
\bibitem{A.G.G}
V.~Apostolov, P.~Gauduchon, G.~Grantcharov,
{\it Bihermitian structures on complex surfaces},
Pro. London Math. Soc. {\bf 79} (1999), 414-428,
Corrigendum: {\bf 92}(2006), 200-202,  MR 2192389, Zbl 1089.53503
\bibitem{B.M}
C.~ Bartocci, E.~ Macr$\grave{\text{\rm i}}$,
{\it Classification of Poisson surfaces},
math.AG/0402338, Commun. Contemp. Math. 7 (2005), no. 1, 89--95. MR 2129789, Zbl 1071.14514
\bibitem{Bo}
Bogomolov,~F.A,
{\it Hamiltonian K\"ahler manifolds}, (English. Russian original), Sov. Math., Dokl. 19, 1462-1465 (1978); translation from Dokl. Akad. Nauk SSSR 243, 1101-1104 (1978). Zbl 0418.53026 
\bibitem{De}
M.~Demazure, {\it Surfaces de del Pezzo II, III, IV, V}, 
 Lecture Notes in Math. 777, pp. 23-69, S\'eminaire sur les Singularit\'es des Surfaces, Palaiseau, France, 1976-1977
\bibitem{DV}
P.~Du Val,
{\it On isolated singularities of surfaces which do not affect the conditions of adjunction. I, II, III}, 
Proc. Cambridge Phi. Soc. 30 (1934), 453-465, 483-491
\bibitem{F.P}
A.~Fujiki and M.~ Pontecorvo,
Bihermitian anti-self-dual structures on Inoue surfaces, 
arXiv:0903.1320, to appear in J.D.G
\bibitem{G.H.R}
S.~J. Gates C.~M. Hull and M. Ro\v cek,
{\it Twisted multiplets and new supersymmetric nonlinear $\sig$ models},
Nuclear Phys. B 248 (1984), 154-186
\bibitem{G.K}
V.~Ginzburg and D.~ Kaledin,
{\it Poisson deformations of symplectic quotient singularities},
Adv. Math. 186 (2004), no. 1, 1--57. 
MR 2065506, Zbl 1062.53074
\bibitem{Go-1}
R.~Goto,
{\it Moduli spaces of topological calibrations,
Calabi-Yau, hyperK\"ahler, G$_2$ and Spin$(7)$ structures},
Internat. J. Math. 15 (2004), no. 3, 211--257. MR 2060789, Zbl 1046.58002
\bibitem{Go0}
R.~Goto,
{\it On deformations of generalized Calabi-Yau, hyperK\"ahler, G$_2$ and Spin$(7)$ structures},
Math.DG/0512211
\bibitem{Go1}
R.~Goto,
{\it Deformations of \complex  and generalized K\"ahler structures},
to appear in J. Differential Geom., 
Math. DG/0705.2495
\bibitem{Go2}
R.~Goto,
{\it Poisson structures and generalized K\"ahler submanifolds},
 J. Math. Soc. Japan 61 (2009), no. 1, 107--132. MR 2272873, Zbl 1160.53014
 \bibitem{Go3}
 R.~Goto,
 {\it Deformations of Generalized K\"ahler structures and Bihermitain structures},
 Math. DG/0910.1651
 \bibitem{Gu1}
M.~Gualtieri,
{\it Generalized complex geometry},
Math.DG/0703298
\bibitem{Hi1}
N.~Hitchin,
{\it Generalized Calabi-Yau manifolds},
Math. DG/0401221, Q. J. Math. 54 (2003), no. 3, 281--308. MR 2013140, Zbl 1076.3201
\bibitem{Hi2}
N.Hitchin,
{\it  Instantons, Poisson structures and generalized K\"ahler geometry},
Comm. Math. Phys. 265 (2006), no. 1, 131--164. MR 2217300, Zbl 1110.53056
\bibitem{Hi3}
N.Hitchin,
{\it Bihermitian metrics on Del Pezzo surfaces},
Math.DG/0608213,
J. Symplectic Geom. 5 (2007), no. 1, 1--8. MR 2371181, Zbl pre05237677
\bibitem{Hu}
D.~Huybrechts,
{\it Generalized Calabi-Yau structures, K3 surfaces and B-fields}, math.AG/0306132,
Internat. J. Math. 16 (2005), no. 1, 13--36. MR 2115675, Zbl 1120.14027
\bibitem{Kob}
P.~Kobak,
{\it Explicit doubly-Hermitian metrics},
Differential Geom. Appl. 10 (1999), no. 2, 179--185, MR 1669453, Zbl 0947.53011
\bibitem{Ko}
K.~Kodaira,
{\it Complex manifolds and deformations of complex structures},
Grundlehren der Mathematischen Wissenschaften [Fundamental Principles of Mathematical Sciences], 283. Springer-Verlag, New York, 1986. x+465 pp. ISBN: 0-387-96188-7 MR 0815922, Zbl 0581.32012
Grundlehren der Mathematischen Wissenschaften, {\bf 283}, springer-Verlag,
(1986)
\bibitem{K.S.I,II}
K.~Kodaira and D.C.~Spencer,
{\it  On deformations of complex, analytic structures I,II},
Ann. of Math. (2) 67 (1958) 328--466. MR 0112154, Zbl 0128.16901
\bibitem{K.S III}
K.~Kodaira and D.C.~Spencer,
{\it On deformations of complex analytic structure, III.
stability theorems for complex structures},
Ann. of Math. (2) 71 (1960) 43--76. MR 0115189, Zbl 0128.16902 
\bibitem{Kos}
Y.~Kosmann-Schwarzbach,
{\it Derived brackets}, 
Lett. Math. Phys. 69 (2004), 61--87. MR 2104437, Zbl 1055.17016 
\bibitem{Ro}
U.~Lindstr\"om, M~. Ro\v cek, Martin, R.~ von Unge and M.~ Zabzine, 
{\it Generalized K\"ahler geometry and manifest $N=(2,2)$ supersymmetric nonlinear sigma-models}, J. High Energy Phys. 2005, no. 7, 067, 18 pp. 
MR2163246.

\bibitem{L.T}
Y.~Lin and S.~Tolman,
{\it Symmetries in generalized K\"ahler geometry},
Comm. Math. Phys. 268 (2006), no. 1, 199--222. MR 2249799, Zbl 1120.53049
, Math. DG/0509069
\bibitem{L.W.P}
Z.-J.~Liu, A.~Weinstein and Ping. Xu, 
{\it Manin triples for Lie bialgebroids},
J. Differential Geom. 45 (1997), no. 3, 547--574. MR 1472888, Zbl 0885.58030
\bibitem{Miya}
K.~Miyajima, 
{\it A note on the Bogomolov-type smoothness on deformations of the regular parts of isolated singularities}, 
Proc. Amer. Math. Soc. 125 (1997), no. 2, 485--492, MR 1363432,
Zbl 0861.32012. 
\bibitem{Na}
Y.~Namikawa,
{\it 
Flops and Poisson deformations of symplectic varieties},
 Publ. Res. Inst. Math. Sci. 44, No. 2, 259-314 (2008),
 MR 2426349, Zbl 1148.14008. 
\bibitem{Pol}
A.~Polishchuk,
{\it Algebraic geometry of Poisson brackets}, Algebraic geometry, 7. J. Math. Sci. (New York) 84 (1997), no. 5, 1413--1444. MR 1465521, Zbl 0995.37057
\bibitem{Sa}
F.~Sakai,
{\it Anti-Kodaira dimension of ruled surfaces},
 Sci. Rep. Saitama Univ. Ser. A 10, no. 2, 1--7. (1982), MR 0662405, Zbl 0496.14022
\bibitem{Sa1}
F.~Sakai, 
{\it Anticanonical models of rational surfaces}, Math. Ann. 269 (1984), no. 3, 389--410. MR 0761313, Zbl 0533.14016
\bibitem{Ti}
G.~Tian,
{\it Smoothness of the universal deformation space of compact Calabi-Yau manifolds and its Petersson-Weil metric}, Mathematical aspects of string theory (San Diego, Calif., 1986), 629--646, Adv. Ser. Math. Phys., 1, World Sci. Publishing, Singapore, 1987. MR 0915841, Zbl 0696.53040
\end{thebibliography}
\end{document}